%% file: FEM_FCT_Solvers.tex
\documentclass{article}
\usepackage[a4paper, total={6.5in, 8.5in}]{geometry}
\usepackage{subcaption}
 \usepackage{color}
\usepackage{graphicx}
\usepackage{amsmath, amssymb, amsthm} 
\usepackage{lmodern} 
\usepackage{microtype}
\usepackage[numbers,square,sort&compress]{natbib}
\usepackage{multirow}
\usepackage{longtable}
\usepackage{url}
\input{new_command.tex}

\title{An Assessment of Solvers for Algebraically Stabilized Discretizations of Convection-Diffusion-Reaction Equations}

\author{Abhinav Jha\thanks{Corresponding Author}  \thanks{\emph{RWTH Aachen University, Applied and Computational Mathematics, Schinkelstra\ss e 2, 52062, Aachen, Germany}, \texttt{jha@acom.rwth-aachen.de}}, Ond\v{r}ej P\'{a}rtl \thanks{\emph{Weierstrass\  Institute for Applied Analysis and Stochastics (WIAS), Mohrenstr.  39, 10117 Berlin, Germany}, \texttt{ondrej.partl@wias-berlin.de}}, Naveed Ahmed\thanks{\emph{Gulf University for Science \& Technology,Block 5, Building 1,Mubarak Al-Abdullah Area, West Mishref Kuwait},  \texttt{ahmed.n@gust.edu.kw}}, and Dmitri Kuzmin \thanks{\emph{Institute of Applied Mathematics (LS III), TU Dortmund University, Vogelpothsweg 87, D-44227 Dortmund, Germany}, \texttt{kuzmin@math.uni-dortmund.de}}}
\date{}

\begin{document}
\maketitle
\begin{abstract}
We consider flux-corrected finite element discretizations of 3D convection-dominated transport problems and assess the computational efficiency of algorithms based on such approximations. The methods under investigation include flux-corrected transport schemes and monolithic limiters. We discretize in space using a continuous Galerkin method and $\mathbb{P}_1$ or $\mathbb{Q}_1$ finite elements. Time integration is performed using the Crank-Nicolson method or an explicit strong stability preserving Runge-Kutta method. Nonlinear systems are solved using a fixed-point iteration method, which requires solution of large linear systems at each iteration or time step. The great variety of options in the choice of discretization methods and solver components  calls for a dedicated comparative study of existing approaches. To perform such a study, we define new 3D test problems for time dependent and stationary convection-diffusion-reaction equations. The results of our numerical
experiments illustrate how the limiting technique,
time discretization and solver impact on the overall performance.
\\\textbf{Classifications:} 65M12, 65M15, 65M60 
\end{abstract}

\section{Introduction}
\input{introduction.tex}

\input{prelim.tex}
\label{sec:preliminaries}

\section{Numerical Studies}
\label{sec:numres}
\input{numres.tex}

\bibliographystyle{abbrv}
\bibliography{FEM_FCT_Solvers}
\end{document}

%% file: new_command.tex
\newcommand{\bfb}{\boldsymbol b}
\newcommand{\bfx}{\boldsymbol x}
\newcommand{\bfn}{\boldsymbol n}
\newcommand{\db}{_\text{D}}
\newcommand{\nb}{_\text{N}}
\newcommand{\bfu}{\boldsymbol u}
\newcommand{\bff}{\boldsymbol F}
\newcommand{\bfr}{\boldsymbol r}
\newcommand{\bfg}{\boldsymbol G}

\newcommand{\mss}{\mathrm{ss}}

\newtheorem{remark}{Remark}

\usepackage{tikz}
\usetikzlibrary{positioning} 
\usetikzlibrary{calc}
\usepackage{pgfplots}
\pgfplotsset{width=10cm,compat=1.9}
\definecolor{light_gray}{gray}{0.75}
\definecolor{lighter_gray}{gray}{0.5}
\colorlet{light_blue}{blue!20}
\definecolor{dark_green}{rgb}{0.0, 0.6, 0.0}
\definecolor{royal_blue}{rgb}{0.0, 0.22, 0.66}
\definecolor{salmon}{rgb}{1.0, 0.55, 0.41}
\definecolor{gold}{rgb}{0.8, 0.63, 0.21}
\definecolor{navy_blue}{rgb}{0.0, 0.0, 0.5}
\definecolor{crimson}{rgb}{0.79, 0.0, 0.09}
\definecolor{amethyst}{rgb}{0.6, 0.4, 0.8}
\definecolor{alizarin}{rgb}{0.82, 0.1, 0.26}
\definecolor{amaranth}{rgb}{0.9, 0.17, 0.31}
\definecolor{azure}{rgb}{0.0, 0.5, 1.0}
\definecolor{canaryyellow}{rgb}{0.82, 0.41, 0.12}
\definecolor{carrotorange}{rgb}{0.8, 0.33, 0.0}
\definecolor{cadmiumgreen}{rgb}{0.0, 0.42, 0.24}
\definecolor{copper}{rgb}{0.72, 0.45, 0.2}
\definecolor{aqua}{rgb}{0.5, 1.0, 0.83}
\definecolor{awesome}{rgb}{1.0, 0.13, 0.32}
\definecolor{candyapplered}{rgb}{1.0, 0.03, 0.0}
\definecolor{caribbeangreen}{rgb}{0.0, 0.8, 0.6}

%% file: introduction.tex
Traditional stabilization techniques for finite element discretizations of convection-diffusion-reaction(CDR) equations
do not ensure the validity of discrete maximum principles
\cite{JS08}. As a consequence, numerical solutions
may attain physically unrealistic values, and simulations
may crash.

In the context of finite volume schemes and
discontinuous Galerkin methods, the relevant inequality
constraints are commonly enforced by using \emph{limiters}
for numerical fluxes or for slopes of piecewise-polynomial
approximations. The first extensions of such schemes to 
continuous finite element approximations \cite{Loh87,Ku02} were based on
generalizations of Zalesak's flux-corrected transport
(FCT) algorithm \cite{Zal79}.

During the last two decades,
many alternatives were developed using the concept of
algebraic flux correction (AFC). The AFC methodology
\cite{KM05,BJK16}  is based on an algebraic splitting
of a high-order target scheme into a
bound-preserving low-order
approximation and an antidiffusive
correction term. The
latter is decomposed into numerical fluxes that are limited to
preserve important properties of the low-order method.

A~theoretical framework for analysis
and design of AFC schemes was developed in \cite{BJK16,Loh19}
and used to construct improved limiter functions in \cite{BJK17,Loh19}.
Further remarkable recent advances in the field include the development
of a monolithic convex limiting strategy for nonlinear
hyperbolic conservation laws and systems \cite{Ku20}.

In this work, we focus on efficient numerical solution of
the nonlinear discrete problems that arise from the AFC discretizations
of time-dependent and stationary problems. The overall
computational cost depends on the type of the limiting
strategy (predictor-corrector vs.~monolithic) and on the time
discretization (explicit vs.~implicit). 
The efficiency of the implicit schemes and steady-state solvers
depends on the convergence rates of the inner and outer iterations \cite{JJ19}.

Unfortunately, these important aspects have received little
attention in the AFC literature so far. We are not aware
of any systematic numerical study focused on the overhead
cost of the flux limiting and the performance--accuracy ratio. However,
the practical use of AFC tools in simulation software for real-life 
applications requires a deeper understanding of such aspects. As
a first step toward that end, we introduce new three-dimensional
test problems and solve them using the AFC schemes proposed in
\cite{BJK17,Ku07,Ku20,JN12}. A direct comparison of CPU times
for different approaches
enables us to identify algorithms that offer the best performance
for a certain class of problems.

This paper has the following structure: The numerical schemes that we study are described in Section \ref{prelim:tcdr} (for the evolutionary problems) and Section \ref{prelim:cdr} (for the stationary problems). Our tests are in Section \ref{sec:numres}.

%% file: prelim.tex
\section{Evolutionary Convection-Diffusion-Reaction Equations}
\label{prelim:tcdr}
We consider the following initial-boundary value problem for a scalar evolutionary convection-diffusion-reaction equation: Find $u:(0,T]\times \Omega \rightarrow \mathbb{R}$ such that
\begin{equation} \label{eq:tcdr}
\begin{array}{rcll}
u_t - \varepsilon \Delta u + \bfb \cdot \nabla u + c u &=& f & \mbox{ in } (0,T] \times \Omega,\\
u  &=& u\db & \mbox{ on } [0,T]\times \Gamma_D,\\
\varepsilon \nabla u \cdot \bfn  &=& g\nb & \mbox{ on } [0,T]\times \Gamma_N,\\
u(0,\bfx) &=& u_0(\bfx ) &\mbox{ } \forall \bfx \in \overline{\Omega}.
\end{array}
\end{equation}
Here $\Omega \subset \mathbb{R}^3$ is a bounded polyhedral domain, $\bfn$ is the outward pointing unit normal to the boundary $\Gamma =\Gamma\db\cup \Gamma\nb$, $\Gamma\db \cap \Gamma\nb=\emptyset$, and \([0,T]\) is a bounded time interval. Furthermore, $\varepsilon$, $0<\varepsilon\ll 1$, is a diffusivity constant, $\bfb = \bfb (t,\bfx)$ denotes a solenoidal velocity field, $c = c (t,\bfx)$ stands for a nonnegative reaction coefficient, and $f=f(t,\bfx)$ represents outer sources of the unknown scalar quantity $u$. On $\Gamma\db$, the Dirichlet boundary conditions ($u\db$) are set, and on $\Gamma\nb$, the Neumann boundary conditions ($g\nb$) are prescribed.

A standard finite element discretization of~\eqref{eq:tcdr} with $\mathbb{P}_1$ or $\mathbb{Q}_1$ elements leads to a system of differential algebraic equations of the form
\begin{equation}\label{eq:high_order_disc}
M_C \dot{\bfu} + A\bfu = \bff,
\end{equation}
where $M_C=\lbrace m_{ij}\rbrace_{i,j=1}^N$ is the consistent mass matrix, $A=\lbrace a_{ij}\rbrace_{i,j=1}^N$ is the stiffness matrix, and $\bff$ is the corresponding right-hand side. 
The length of the vectors is denoted by $N$, which corresponds to the number of degrees of freedom. The matrix entries are given by
\begin{align}
m_{ij} & = \left( \varphi_j, \varphi_i\right), \\
\label{eq:a_ij_initial_assembling}
a_{ij} & = \varepsilon \left(\nabla \varphi_j, \nabla \varphi_i\right) + \left( \bfb \cdot \nabla \varphi_j,\varphi_i\right)+\left(c\varphi_j,\varphi_i\right),
\end{align}
where $(\cdot,\cdot)$ denotes the standard inner product in $L^2 (\Omega)$, and $\lbrace \varphi_i \rbrace_{i=1}^N$ is the standard finite element basis.

The first step in the AFC methodology is to modify~\eqref{eq:high_order_disc} so that we obtain an M-matrix $\mathbb{A}$ instead of $A$. For this purpose, we define the lumped mass matrix $M_L$ and an artificial diffusion matrix $D$ as follows:
\begin{align}
M_L & = \mathrm{diag}(m_i), \quad m_i=\sum_{j=1}^Nm_{ij},\nonumber \\
D & = \lbrace d_{ij}\rbrace_{i,j=1}^N, \quad d_{ij}=-\max \lbrace a_{ij}, 0, a_{ji}\rbrace \ \text{ for } i\neq j,\ \ \ d_{ii}=-\sum_{j=1,j\neq i}^Nd_{ij}.
\end{align}
Replacing $M_C$ by $M_L$ and $A$ by $\mathbb{A}=A+D$ in~\eqref{eq:high_order_disc}, we obtain
\begin{equation}\label{eq:low_order_disc}
M_L \dot{\bfu} + \mathbb{A}\bfu=\bff.
\end{equation}
The temporal discretization of this equation yields a low-order scheme that is bound preserving, but overly diffusive.

To reduce this excessive diffusivity, we add an antidiffusive correction term $\bff^*$ on the right-hand side of \eqref{eq:low_order_disc} to get
\begin{eqnarray}\label{eq:flux_represnt_001}
M_L \dot{\bfu} + \mathbb{A}\bfu = \bff+\bff^* .
\end{eqnarray}
To define $\bff^*$, we first consider the residual difference $\bfr$ obtained by  subtracting~\eqref{eq:high_order_disc} from~\eqref{eq:low_order_disc}:
\begin{eqnarray}\label{eq:flux_000}
\bfr & = & \left( M_L-M_C\right) \dot{\bfu} + D\bfu.
\end{eqnarray}
Next, we decompose each component of $\bfr$ as
\begin{equation}\label{eq:r_decomposition}
r_i = \sum_{j=1,j\neq i} r_{ij}, \ \text{ where } \ r_{ij} = m_{ij} ( \dot{u}_i - \dot{u}_j ) + d_{ij} ( u_j - u_i ).
\end{equation}
Using this decomposition, we set
\begin{eqnarray}\label{eq:definition_of_f_i}
F_i^* & = & \sum_{j=1,j\neq i}^N\alpha_{ij}r_{ij},
\end{eqnarray}
where $\{\alpha_{ij}\}_{ij=1}^N \subset [0, 1]$ are solution-dependent correction factors; algorithms for calculating them are called limiters. For $\alpha_{ij}=1$, we revert to the standard Galerkin formulation, whereas setting $\alpha_{ij}=0$ corresponds to the over-diffusive scheme. Note that the computation of $r_{ij}$ is no longer required when $i$ is a Dirichlet node.

Various definitions of $\alpha_{ij}$ have been proposed in the literature (see \cite{Zal79, Ku12, BJK17, Ku20}). Some flux limiters are defined at the semi-discrete
level and applicable to steady-state problems as well. In other approaches, the fluxes
$r_{ij}$ and the correction factors $\alpha_{ij}$ are derived for a particular
time-stepping method. In the next sections, we introduce the time
discretizations and limiters used in our numerical studies.

\subsection{Flux-Corrected Transport Algorithms}

We begin with algorithms that use (generalizations of) Zalesak's FCT limiter
\cite{Zal79} to calculate the correction factors $\alpha_{ij}$.
The anti\-diffusive fluxes $r_{ij}$ corresponding
to specific time integrators are defined below.

\subsubsection{Crank-Nicolson Scheme}

The Crank-Nicolson (CN) time discretization of the semi-discrete problem
\eqref{eq:flux_represnt_001} yields the nonlinear system
\begin{eqnarray}\label{CNZal}
  \left[\frac{1}{\Delta t}M_L +\frac{1}2\mathbb{A}\right]\bfu^n =
 \left[\frac{1}{\Delta t}M_L -\frac{1}2\mathbb{A}\right]\bfu^{n-1}
 + \frac{1}{2}\bff^n + \frac{1}{2}\bff^{n-1} + \bff^*(\bfu^n,\bfu^{n-1}),
\end{eqnarray}
where $\Delta t$ is the time-step length, the superscripts denote the
time levels, and the correction term
$\bff^*$ is assembled from limited counterparts
$\alpha_{ij}r_{ij}$ of the antidiffusive fluxes \cite{KM05,Ku09}
\begin{equation}\label{eq:flux_002}
r_{ij}= \frac{m_{ij}}{\Delta t} \left[ u_i^n-u_i^{n-1}-\left( u_j^n-u_j^{n-1}\right) \right]
+\frac{d_{ij}}{2}\left[ u_j^n+u_j^{n-1}-\left( u_i^n+u_i^{n-1}\right)\right].
\end{equation}
The CN-Galerkin scheme is a nondissipative high-order method, which
tends to generate small ripples within the local bounds of the
limiting procedure \cite{KM05}. This behavior can often be cured by using a high-order
linear stabilization or prelimiting \cite{Ku09}. In this work, we
 prelimit $r_{ij}$ as follows:
\begin{equation}\label{eq:prelimiting}
r_{ij} = \text{minmod} \left( r_{ij}, \, L \, d_{ij} ( u_j - u_i ) \right),
\end{equation}
where
\begin{equation*}
\text{minmod} ( a, \, b ) =
\begin{cases} 0 \ \ \text{if} \ \ ab < 0 , \\
\min \left\lbrace a,b \right\rbrace \ \ \text{if} \ \  a>0 \land b>0 , \\
\max \left\lbrace a,b \right\rbrace \ \ \text{if} \ \  a<0 \land b<0 ,
\end{cases}
\end{equation*}
and $L = 2$ is a Lipschitz constant based on the analysis in \cite{BJK17}.

In addition to the fluxes $r_{ij}$, Zalesak's limiter (as presented
in Section \ref{prelim:zalesak,tcdr} below) requires a bound-preserving
intermediate solution $\tilde{\bfu}$ of low order. For the CN version, it is defined by
\begin{equation}\label{eq:forward_euler}
\tilde{\bfu} = \bfu^{n-1}-\frac{\Delta t}{2}M_L^{-1}\left( \mathbb{A}\bfu^{n-1}-\bff^{n-1}\right),
\end{equation}
which can be viewed as the solution of~\eqref{eq:low_order_disc} at time $t_{n-1/2}$ computed using the forward Euler scheme with time step $\Delta t/2$. Note that $\tilde{\bfu}$ should be constrained to satisfy the
Dirichlet boundary conditions for nodes belonging to $\Gamma_D$.

\begin{remark}
  As shown in \cite{KM05,Loh19}, the explicit predictor $\tilde{\bfu}$ is bound
  preserving under a CFL-like condition. 
\end{remark}

We test two implementations of the CN-FCT algorithm: The first one
solves the system of equations for $\bfu^n$ using a fixed point iteration
method, which means that a sparse linear system needs to be solved at each step.
A brief overview of the iterative procedure is given in Section~\ref{sec:numres}. In what follows, we refer to this scheme as \textit{nonlinear Zal+CN}.

The second algorithm is the linearized CN scheme proposed in \cite{JN12} that utilizes
$\tilde{\bfu}$ defined by~\eqref{eq:forward_euler} to approximate $\bfu^{n}$ by
$2 \tilde{\bfu} - \bfu^{n-1}$ in formula \eqref{eq:flux_002}. Hence, it replaces
\eqref{CNZal} by a linear system for $\bfu^n$. In the following sections, we refer to this scheme as \textit{linear Zal+CN}.

\subsubsection{Second-Order SSP Scheme}

As an alternative to the implicit CN scheme, we consider the second-order
explicit strong stability preserving (SSP) time
integrator commonly known as Heun's method. For a linear or
nonlinear system of the form

\begin{equation}\label{SSP-DAE}
\dot{\bfu}(t)=\bfg \left(\bfu (t), t\right),
\end{equation}
the numerical solution at time $t_n$ is given by
\begin{equation*}
  \bfu^n = \bfu^{n,0} + \frac{\Delta t}2\left[
    \bfg \left(\bfu^{n,0} , t^{n,0} \right)+
    \bfg \left(\bfu^{n,1}, t^{n,1} \right)\right]
  =\frac{\bfu^{n,0}+\bfu^{n,2}}2,
\end{equation*}
where
\begin{align}
  &t^{n,0} = t^{n-1},\quad t^{n,1} = t^n , \\
  &\bfu^{n,0}=\bfu^{n-1},\quad
  \bfu^{n,s}=\bfu^{n,s-1} + \Delta t  \bfg \left(\bfu^{n,s-1} , t^{n,s-1} \right),\qquad
  s=1,2.
  \end{align}
The application of this method to \eqref{eq:flux_represnt_001} requires
two explicit Euler updates of the form
\begin{equation}\label{FE-SSP}
 \bfu^{\rm new}  =\bfu+\Delta tM_L^{-1} \left(
 \bff+\bff^*(\bfu)-\mathbb{A}\bfu \right).
\end{equation}

At each stage, the flux limiting is performed using Zalesak's algorithm
with the low-order predictor

\begin{equation}
\tilde\bfu=\bfu+\Delta tM_L^{-1} \left(
 \bff-\mathbb{A}\bfu \right)
\end{equation}
and the antidiffusive fluxes
\begin{equation}\label{eq:defn_r_ij_SSP}
r_{ij}=m_{ij}\left( \dot{u}^L_i-\dot{u}^L_j\right) + d_{ij}\left( u_j - u_i\right),
\end{equation}
where
\begin{equation}\label{eq:lower_time_approx}
\dot{\bfu}^L=M_L^{-1}\left( \bff - \mathbb{A} \bfu \right)
\end{equation}
stands for the low-order approximation given by~\eqref{eq:low_order_disc}.
As shown in \cite{Ku09}, the fluxes defined by \eqref{eq:defn_r_ij_SSP}
do not require prelimiting because the use of the low-order time
derivatives introduces high-order linear stabilization.

\subsubsection{Zalesak's Limiter}
\label{prelim:zalesak,tcdr}

Given a low-order predictor $\tilde\bfu$ and an array of antidiffusive
fluxes $r_{ij}$, our FCT schemes
use Zalesak's limiter \cite{Zal79}
to calculate the correction factors $\alpha_{ij}$ as follows:
\begin{enumerate}
    \item Compute 
    \begin{equation*}
       P_i^+ = \sum_{j=1,j\neq i}^N \max \lbrace r_{ij},0\rbrace, \qquad P_i^- = \sum_{j=1,j\neq i}^N \min \lbrace r_{ij},0\rbrace.
    \end{equation*}
    \item Compute 
    \begin{eqnarray}
       Q_i^+ & = &\max\left\{ 0,\max_{j=1,\ldots,N,j\neq i}(\tilde{u}_j-\tilde{u}_i) \right\},\nonumber \\
       Q_i^- & = &\min\left\{ 0,\min_{j=1,\ldots,N,j\neq i}(\tilde{u}_j-\tilde{u}_i) \right\}.\nonumber
    \end{eqnarray}
    \item Compute
    \begin{equation*}
       R_i^+ = \min\left\{ 1,\frac{m_iQ_i^+}{\Delta t P_i^+} \right\}, \qquad R_i^- = \min\left\{ 1,\frac{m_iQ_i^-}{\Delta t P_i^-} \right\}.
    \end{equation*}
    If \(P_i^+\) or \(P_i^-\) is zero, we set \(R_i^+=1\) or \(R_i^-=1\), respectively. We also set $R_i^+ = R_i^- = 1$ if $i$ is a Dirichlet node. 
    \item Compute
    \[ \alpha_{ij} = \begin{cases} 
        \min\{R_i^+,R_j^-\} & \mbox{ if } r_{ij}>0,\\
        \min\{R_i^-,R_j^+\} & \mbox{ otherwise}.
      \end{cases}
   \]
\end{enumerate}
The CN and SSP versions of Zalesak's FCT scheme
differ in the definition of $r_{ij}$ and
$\tilde\bfu$.

\subsection{Monolithic Convex (MC) Limiter}
\label{sec:MCL}

A potential drawback of FCT-like approaches is their dependence
on the particular time-stepping method. As an alternative that
is also applicable to stationary problems, we consider the
monolithic convex (MC) limiting algorithm proposed in \cite{Ku20}.
Flux limiters of this kind use definition \eqref{eq:defn_r_ij_SSP}
of $r_{ij}$ for the semi-discrete scheme \eqref{eq:flux_represnt_001}. The
limited antidiffusive term of the CN scheme \eqref{CNZal} is given by
$$
F_i^*(\bfu ^n,\bfu^{n-1})
=\frac{1}{2}\sum_{j=1,j\neq i}^N \left( \alpha_{ij}(\bfu ^n)
r_{ij}(\bfu^n)+\alpha_{ij}(\bfu^{n-1})r_{ij}(\bfu^{n-1})\right),
$$
while the SSP version performs forward Euler updates \eqref{FE-SSP} using
$$
F_i^*(\bfu)=\sum_{j=1,j\neq i}^N\alpha_{ij}(\bfu)r_{ij}(\bfu).
$$
The corresponding nonlinear space discretizations are of the form
\eqref{SSP-DAE} and reduce to
$\bfg \left(\bfu\right)=0$ at steady state. At each
fixed-point iteration or Runge-Kutta stage, the limited fluxes
$\alpha_{ij} r_{ij} = r^*_{ij}$ are defined by
\begin{equation}\label{MClimiter}
r^*_{ij} =
	\begin{cases} 
	\min\left\lbrace r_{ij}, \min \left\lbrace 2d_{ij} \left(\bar{u}_{ij}-u_i^{\max}\right), 2d_{ij} \left(u_j^{\min}-\bar{u}_{ji}\right) \right\rbrace \right\rbrace \text{ if } r_{ij}>0, \\
	\max\left\lbrace r_{ij}, \max \left\lbrace 2d_{ij} \left(\bar{u}_{ij}-u_i^{\min}\right), 2d_{ij} \left(u_j^{\max}-\bar{u}_{ji}\right)\right\rbrace \right\rbrace \text{ otherwise},
	\end{cases}
\end{equation}
where $\bar u_{ij}$ are intermediate states defined by
$$2d_{ij} \bar{u}_{ij}  =  d_{ij} ( u_i+u_j ) + a_{ij} (u_j-u_i)$$
and
\begin{align}
\label{eq:u_max_min}
u_i^{\max} = \max_{j\in N_i}u_j,\ \ u_i^{\min} = \min_{j\in N_i}u_j.
\end{align}
In the last formula, 
$N_i=\{j\in\{1,\ldots,N\}\,:\,m_{ij}\ne 0\}$ is the integer set
containing the indices of node $i$ and its nearest neighbors.
Note that the definition of $u_i^{\max}$ and $u_i^{\min}$ can be
changed to ensure linearity preservation 
\cite[Section 6.1]{Ku20}.

\begin{remark}\label{rem:assembling_mc}
  Note that the MC limiter defined by \eqref{MClimiter}
  was designed for hyperbolic conservation laws. When applying this limiter to convection-diffusion-reaction (CDR) equations, we perform algebraic flux correction for the semi-discrete problem corresponding to $\varepsilon=0,\ c=0$ and add the unlimited discretization of $\varepsilon \Delta u-cu$ on the right-hand side of the resulting system.
  
However, a proper extension of the MC limiter to CDR problems should include the diffusive and reactive terms in a manner that ensures preservation of local bounds. The development of such extensions is beyond the scope of the present work which is mainly focused on solver aspects.
\end{remark}

\section{Stationary Convection-Diffusion-Reaction Equations}
\label{prelim:cdr}
In this section, we present AFC schemes for the stationary counterpart of~\eqref{eq:tcdr}, the boundary value problem
\begin{equation}\label{eq:cdr}
\begin{array}{rcll}
- \varepsilon \Delta u + \bfb \cdot \nabla u + c u &=& f & \mbox{ in } \Omega,\\
u  &=& u\db & \mbox{ on } \Gamma_D,\\
\varepsilon \nabla u \cdot \bfn  &=& g\nb & \mbox{ on } \Gamma_N.
\end{array}
\end{equation}

Flux-limited discretizations of such problems always lead to a nonlinear system of equations. We solve these systems via the fixed point iteration studied in \cite{JJ18, JJ19} and outlined at the beginning of Section \ref{sec:numres}. This iterative solver requires the solution of a system of linear
equations at each step.

When deriving our numerical schemes for~\eqref{eq:cdr}, we use the same procedure as in Section \ref{prelim:tcdr}, but the time derivatives vanish. Hence, the antidiffusive term $\bfr$ reduces to $\bfr^{\mss}$, where
\begin{equation}\label{eq:r_ss_decomposition}
r_i^{\mss} = \sum_{j=1,j\neq i} r_{ij}^{\mss} = \sum_{j=1,j\neq i} d_{ij}(u_j-u_i).
\end{equation}

\subsection{Monolithic Convex (MC) Limiter}

The MC limiter presented in Section \ref{sec:MCL} is directly
applicable to stationary problems. At each fixed-point iteration, the antidiffusive fluxes
$r_{ij}^{\mss}=d_{ij}(u_j-u_i)$ are limited
using formula \eqref{MClimiter}. The validity of a discrete maximum
principle for the converged steady-state solution was shown
in \cite[Theorem A.3]{Ku20}.

\subsection{Monolithic Upwind (MU) Limiter}

This limiter was proposed in \cite{Ku07} and analyzed in \cite{BJK16}.
Using the notation
\begin{equation}\label{eq:rijss+-}
r_{ij}^{\mss,+} = \max \left\lbrace r^{\mss}_{ij}, 0 \right\rbrace, \ \ \ r_{ij}^{\mss,-} = \min \left\lbrace r^{\mss}_{ij}, 0 \right\rbrace,
\end{equation}
the correction factors $\alpha_{ij}$ are computed as follows: 
\begin{enumerate}
    \item Compute 
    \begin{equation*}
       P_i^+ = \sum_{j=1, \, a_{ji}\leq \, a_{ij}}^N r_{ij}^{\mss,+}, \qquad P_i^- = \sum_{j=1, \,a_{ji}\leq \, a_{ij}}^N r_{ij}^{\mss,-}.
    \end{equation*}
    \item Compute 
       \begin{equation*}
       Q_i^+ = -\sum_{j=1}^N r_{ij}^{\mss,-}, \qquad Q_i^- = -\sum_{j=1}^N r_{ij}^{\mss,+}.
    \end{equation*}
    \item Compute
    \begin{equation*}
       R_i^+ = \min\left\{ 1,\frac{Q_i^+}{P_i^+} \right\}, \qquad R_i^- = \min\left\{ 1,\frac{Q_i^-}{P_i^-} \right\}.
    \end{equation*}
    If $i$ is a Dirichlet node or if $P_i^+$ or $P_i^-$ is zero, the corresponding $R_i^+$ or $R_i^-$ is set to 1.
    \item For all $i,j$ such that $a_{ji}\leq a_{ij}$, set
    \begin{align}\label{eq:alphas_symmetry_condition}
    \alpha_{ij} &= \begin{cases} 
        R_i^+ & \mbox{ if } r_{ij}^{ss}>0,\\
        1 & \mbox{ if } r_{ij}^{ss}=0,\\
        R_i^- & \mbox{ if } r_{ij}^{ss}<0,
      \end{cases}\qquad \alpha_{ji} = \alpha_{ij} .
    \end{align}
\end{enumerate}
Similarly to the MC limiter, this algorithm was designed
for the hyperbolic case. It exploits the skew symmetry
of the discrete convection operator and requires a careful
extension to transport problems with diffusion and/or reaction.

\subsection{Linearity Preserving (LP) Limiter}
This limiter, which makes the AFC scheme linearity preserving, was introduced in \cite{BJK17}. It is custom-made for $\mathbb{P}_1$ elements.
The correction factors $\alpha_{ij}$ are computed as follows:
\begin{enumerate}
    \item Compute 
    \begin{equation*}
       P_i^+ = \sum_{j=1,j\ne i}^N r_{ij}^{\mss,+}, \qquad P_i^- = \sum_{j=1,j\ne i}^N r_{ij}^{\mss,-},
    \end{equation*}
    where $r_{ij}^{\mss,+}$ and $r_{ij}^{\mss,-}$ are given by \eqref{eq:rijss+-}.
    \item Compute 
       \begin{equation*}
       Q_i^+ = q_i(u_i^{\max}-u_i), \qquad Q_i^- = q_i(u_i^{\min}-u_i),
    \end{equation*}
    where the bounds $u_i^{\max}$ and $u_i^{\min}$ are defined as in \eqref{eq:u_max_min}, and
    $$
    q_i=-\sum_{j\in N_i}\gamma_i d_{ij}
    $$
    for a positive constant
    $\gamma_i$ depending only on the shape of the spatial grid in the nearest vicinity of the node $i$. We define $\gamma_i$ as in \cite[Rem.~6.2]{BJK17}.
    \item Compute
    \begin{equation*}
       R_i^+ = \min\left\{ 1,\frac{Q_i^+}{P_i^+} \right\}, \qquad R_i^- = \min\left\{ 1,\frac{Q_i^-}{P_i^-} \right\}.
    \end{equation*}
If $i$ is a Dirichlet node or if $P_i^+$ or $P_i^-$ is zero, the corresponding $R_i^+$ or $R_i^-$ is set to 1.
    \item For all $i,j$ define
    \[ \overline{\alpha}_{ij} = \begin{cases} 
        R_i^+ & \mbox{ if } r_{ij}^{ss}>0,\\
        1 & \mbox{ if } r_{ij}^{ss}=0, \\
        R_i^- & \mbox{ if } r_{ij}^{ss}<0.
      \end{cases}
   \]
     For each combination of non-Dirichlet nodes $i$ and $j$, set
   \[ 
   \alpha_{ij} = \min \left\lbrace\overline{\alpha}_{ij}, \overline{\alpha}_{ji}\right\rbrace .
   \]
   For each combination of a Dirichlet node $j$ and a non-Dirichlet node $i$, set
   $$
   \alpha_{ij}=\overline{\alpha}_{ij}.
   $$
\end{enumerate}

The boundary conditions are taken into account by setting
$a_{ij} = 0$ for each combination of a non-Dirichlet
node $i$ and a Dirichlet node $j$.

%% file: numres.tex
In this section, we perform numerical studies for stationary and evolutionary convection-diffusion-reaction problems in 3D domains to evaluate the accuracy and efficiency of the above AFC methods.

The nonlinear systems corresponding to flux-corrected schemes from Sections \ref{prelim:tcdr} and \ref{prelim:cdr} were solved using a fixed point method with dynamic damping. We give a brief overview of the scheme for the stationary problem. The same solution strategy is used for the evolutionary problem. We refer the reader to \cite{JJ19} for a detailed explanation.

The matrix form of the nonlinear schemes under investigation is given by
$$
\mathbb{A}\bfu = \bff+ \bff^*(\bfu).
$$
We solve such nonlinear systems using fixed-point iterations of the form 
$$
\mathbb{A}\bfu^{\nu+1}=\bff +\omega\bff^*(\bfu^{\nu}),
$$
where $\nu$ is the $\nu$-th iteration step, and $\omega$ is a dynamic damping parameter. Here, the matrix $\mathbb{A}$ is a constant $M$-matrix, and thus must be factorized once, and the factorization can be reused in the iteration loop. For the CN discretization of a time-dependent problem, the iteration matrix is $((\Delta t)^{-1}M_L +\frac12 \mathbb{A})$. This is also a constant $M$-matrix, and hence can be reused in the iteration process.

To evaluate the efficiency of different limiters, we compare the computation times using the following parameters:
\begin{enumerate}
\item \label{enum:solvers} \textbf{Choice of solvers:} 
We have tested several direct and iterative solvers for linear systems available in the Portable library, Extensible Toolkit for Scientific Computation Toolkit for Advanced Optimization (\textsc{PETSc}), Release version 3.14.3  \cite{BAA19, PETSc_1, PETSc_2}. The solvers, along with the \textsc{PETSc} arguments and abbreviations used in the simulations, are listed in Tables~\ref{tab:solvers_iterative} and \ref{tab:solvers_direct}. They are all used with the default settings, except that we set $\mathtt{rtol} = 0$, and we use various values of $\mathtt{atol}$.
This means that the stopping criterion for the linear solvers was
\begin{equation*}
R < \mathtt{atol},
\end{equation*}
where $R$ is the Euclidean norm of the residual. We also specified various maximum numbers of iterations.

\begin{table}[t!]
\centering
\caption{Iterative solvers and the corresponding \textsc{PETSc} arguments.}
\begin{tabular}{l l l c}
\textbf{Name} & \textbf{\textsc{PETSc} arguments} & \textbf{Abbreviation} & \textbf{Marker}\\
 & (\texttt{-ksp\_type}) & & \\ \hline
flexible GMRES & \texttt{fgmres} & FGMRES & $\bullet$\\
loose GMRES & \texttt{lgmres} & LGMRES & $\blacksquare$\\
stabilized BiCG & \texttt{bcgs} & BCGS & $\blacktriangle$\\
\end{tabular}
\label{tab:solvers_iterative}
\end{table}

\begin{table}[t!]
\centering
\caption{Direct solvers and the corresponding \textsc{PETSc} arguments.}
\begin{tabular}{l l l}
\textbf{Name} & \textbf{\textsc{PETSc} arguments} & \textbf{Abbreviation}\\
 & (\texttt{-pc\_factor\_mat\_solver\_type}) &  \\ \hline
LU factorization & \texttt{mumps} & LU \\
UMFPACK & \texttt{umfpack} & UMFPACK \\
\end{tabular}
\label{tab:solvers_direct}
\end{table}

\item \textbf{Choice of preconditioners}: The iterative solvers of Table~\ref{tab:solvers_iterative} are used in combination with various preconditioners from \textsc{PETSc}, see the list of preconditioners along with the corresponding \textsc{PETSc} arguments and abbreviations in Table~\ref{tab:pc}.

\begin{table}[t!]
\centering
\caption{List of preconditioners and \textsc{PETSc} arguments.}
\begin{tabular}{l l l c}
\textbf{Name} & \textbf{\textsc{PETSc} arguments} & \textbf{Abbreviation} & \textbf{Color}\\
&  (\texttt{-pc\_type}) & & \\\hline
Jacobi & \texttt{jacobi} & Jac & \textcolor{black}{gold} \\
point block Jacobi & \texttt{bjacobi} & BJac & \textcolor{black}{red}\\
successive over relaxation & \texttt{sor} & SOR & \textcolor{black}{green}\\
additive Schwarz method & \texttt{asm} & ASM & \textcolor{black}{blue}\\
multigrid & \texttt{mg} & MG & --- \\
\end{tabular}
\label{tab:pc}
\end{table}

\item \label{enum:parallelization} \textbf{Effect of parallelization}: The numerical simulations are performed both sequentially and in parallel and compared in terms of the resulting computing times. The parallel calculations were performed using the \textit{Open Run-Time Environment} (\textsc{OpenRTE}), version 1.10.7.0.5e373bf1fd \cite{castain05:_open_rte}, which originated from the Open MPI project.

In our representation, the number of processors is denoted by NP. We consider NP = 4, 8, 16, and 32. The sequential case is denoted by NP = 1.

Note that the definitions of the preconditioners BJac, SOR, and ASM depend on the number of processors by default (e.g., there is one block for each processor in BJac). We considered this dependence to be so natural that we kept it when we changed the number of processors. So, strictly speaking, if we changed the number of processors, we were no longer using the same solver.
\end{enumerate}

To check the accuracy of the schemes, we compare the $L^1$ and $L^2$ norms (denoted by $\left\Arrowvert \cdot \right\Arrowvert_1$ and $\left\Arrowvert \cdot \right\Arrowvert_2$) of the error $(u - u_h )$, the plots of the solution, and the time evolution of the solution for a given point.

In the presentation of the numerical results, the following is always the same:
\begin{itemize}
\item The spatial grid is always generated by successive uniform refinement starting from an initial grid. The number of refinement steps performed is called the \textit{(refinement) level}.
\item Unless otherwise specified, the termination criterion for the nonlinear solver is
\begin{equation}
\label{eq:tol_formula}
R < \sqrt{\# \text{DOF}} \cdot \mathtt{tol},
\end{equation}
where $\#\mathrm{DOF}$ is the number of degrees of freedom, and $\mathtt{tol}$ is a positive number close to zero.
\item We always use the same pattern when representing the computing times in figures (see, e.g., Figure~\ref{fig:bjk17}, page~\pageref{fig:bjk17}): The lines for BJac are always red, the lines for BCGS are always dotted and so on. Tables~\ref{tab:solvers_iterative} and \ref{tab:pc} contain the colors and the indices used in this work to denote the solvers and
preconditioners.
\end{itemize}

The simulations were performed using the in-house code \textsc{ParMooN} \cite{WB16} on the computer HPE Synergy 660 Gen10 with 2 Xeon eighteen-core processors, 3000 MHz and 768 GB RAM.
\subsection{Time-dependent problems}\label{ex:tcdr}

\subsubsection{Rotating shapes}\label{ex:rotating_shapes}
\input{rotating_shapes}
\clearpage

\subsubsection{Concentration of species}\label{ex:concentration_species2}
\input{concentration_species2}
\clearpage

\subsection{Stationary problems}\label{ex:cdr}
\subsubsection{Example with non-constant convection}\label{ex:non_constant_convection2}
\input{example_non_constant_convection2.tex}

\clearpage

\subsubsection{Circular convection}\label{ex:circular_convection}
\input{circular_convection.tex}
\clearpage

\section{Conclusions}

The numerical studies of AFC schemes in this work indicate that
the costs of the following procedures must be taken into account when
designing high-performance algorithms: (i) calculation of the correction factors for limited
antidiffusive fluxes, (ii) matrix/residual assembly, and (iii)
iterative solution.

The need for high accuracy and efficiency becomes particularly
pronounced in applications of AFC to 3D problems. As a rule,
schemes equipped with more diffusive limiters converge
faster than more accurate approaches. Thus a fair comparison
of different algorithms should be based on CPU times that
are required to attain a certain level of accuracy.

An interesting observation is that implicit schemes can
outperform their explicit counterparts even in numerical
simulations of transient processes using small time steps.
We hope that our comparison of different limiters and linear
algebra tools provides useful insights for finding combinations
of methods that offer the best overall performance in terms
of accuracy and efficiency.

%% file: rotating_shapes.tex
Our first example is inspired by the well-known two-dimensional transport problem with rotating bodies \cite{Le96}. As far as we know, this example has not been treated in the literature before. This problem aims to investigate the accuracy of the newly introduced MC limiter and to compare the results with the Zalesak limiter.

\subsubsection*{Description of Problem}

\begin{figure}[t!]
\centering
\caption{Example~\ref{ex:rotating_shapes}: Initial condition.}
\includegraphics[width=0.5\textwidth, trim = 1cm 2cm 5cm 4cm, clip = true]
{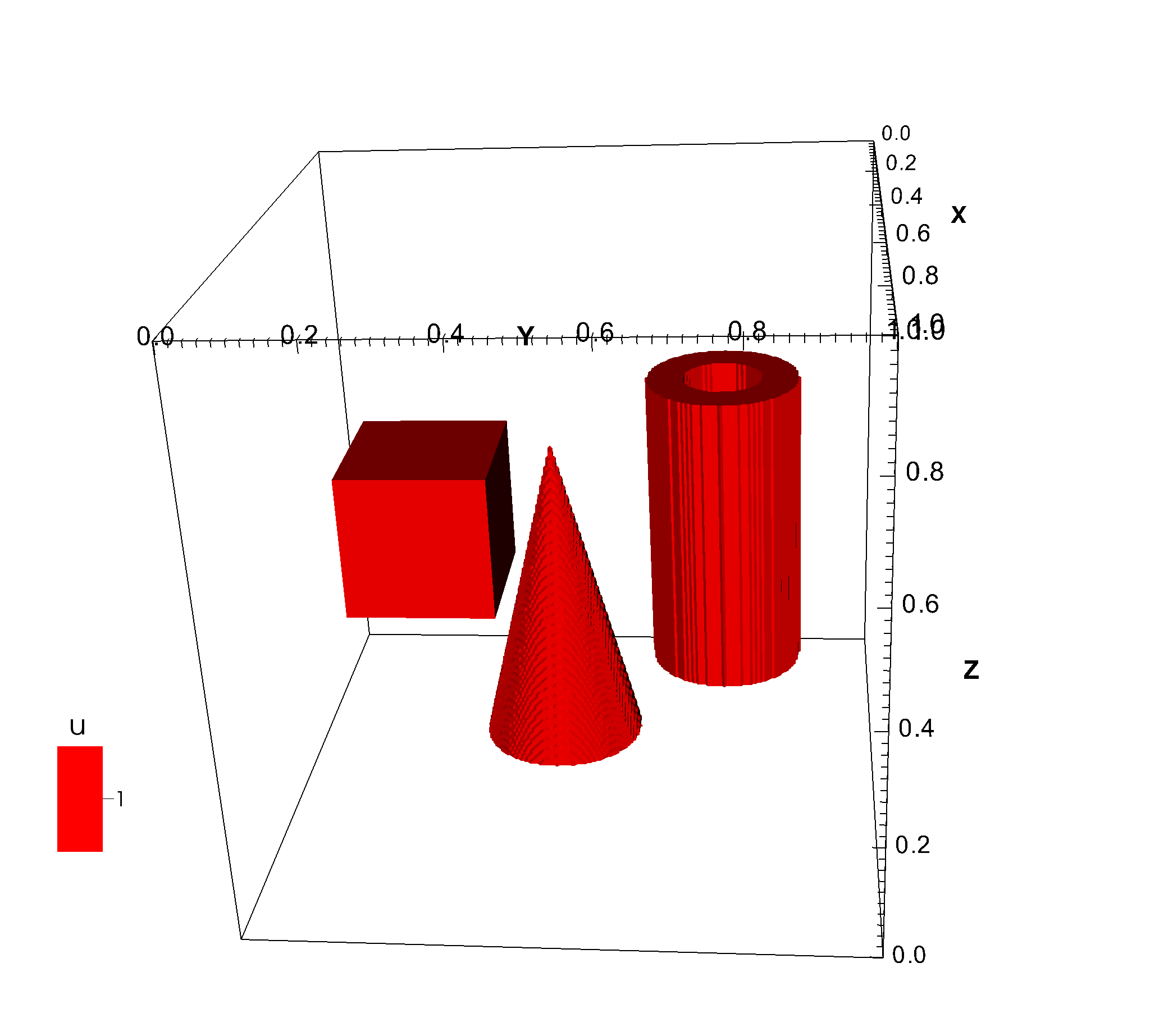}
\label{fig:ini_rotating_shapes}
\end{figure}

We consider~\eqref{eq:tcdr} with $\Omega = (0,1)^3$, $\varepsilon = 0$, $\bfb = (0.5 - y, x-0.5, 0)^\text{T}$, $c = 0$, $f = 0$, $T=2\pi$, $\Gamma\db = \Gamma$, $\Gamma\nb = \emptyset$, $g\nb=0$, and $u\db = 0$. The initial condition $u_0$ is depicted in Figure~\ref{fig:ini_rotating_shapes}, i.e., $u_0$ takes the value 1 in the volumes enclosed by the red surfaces and zero everywhere else. These volumes have the following properties:
\begin{itemize}
    \item cube: edges of length 0.25 parallel to coordinate axes, center at $(0.5, 0.25, 0.5)^T$; 
    \item cone: height 0.5 and bottom surface of radius $0.125$ and center at $(0.75, 0.5, 0.25)^T$ parallel to plane $z=0$;
    \item hollow cylinder: top and bottom surfaces parallel to plane $z=0$, inner radius 0.0625, outer radius 0.125, height 0.5, bottom surface with center at $(0.5, 0.75, 0.25)^T$.
\end{itemize}

The velocity field $\bfb$ makes the spatial shapes rotate around the axis $x = y = 0.5.$ One full rotation takes $t=2\pi$. Since $\varepsilon$ and $f$ both are zero, one should obtain a solution similar to the initial solution after one full rotation.    

Numerical simulations were performed with the uniform cube mesh with the edge length $2^{-8}$ and a fixed time step length $\Delta t = 10^{-3}$. The implicit schemes used LGMRES with SOR as the preconditioner.

Finally, we set $\mathtt{atol} = 10^{-25}$ for the \textsc{PETSc} solver. The nonlinear solver stopped if $R < 10^{-20}$ or after 50 iterations. The reason for this maximum number of iterations is explained below.

\subsubsection*{Discussion of Results}

\begin{figure}[t!]
\centering
\caption{Example~\ref{ex:rotating_shapes}: $u_0(\bfx)$ and the resulting $u(T, \bfx)$ for $\bfx$ in the plane given by $z=0.5$.}
\begin{subfigure}{0.45\textwidth}
\caption{Initial condition.}
\includegraphics[width=\textwidth, trim = 2.5cm 6cm 4cm 2cm, clip = true]
{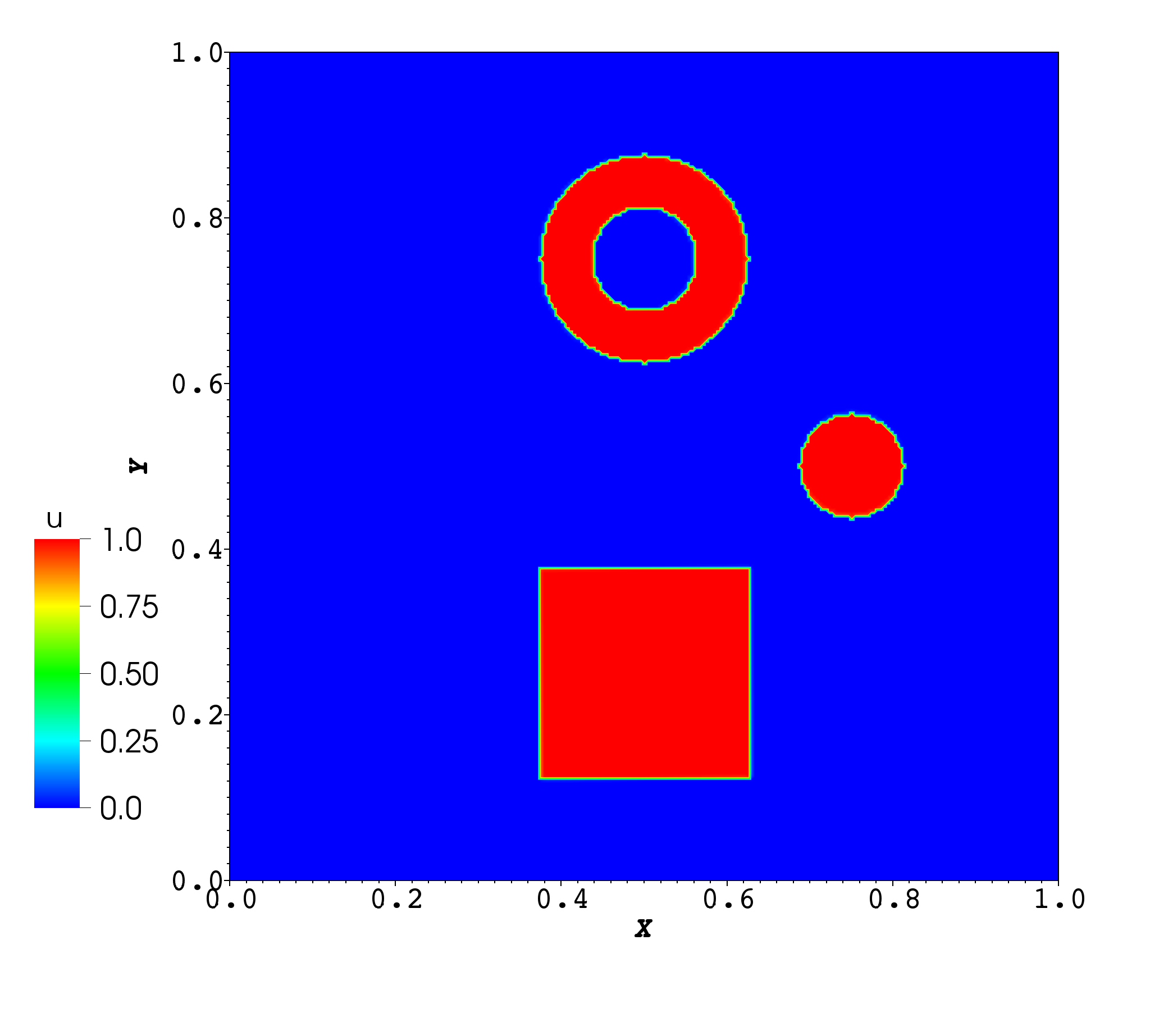}
\label{fig:ini_rotating_shapes_slice}
\end{subfigure}
\begin{subfigure}{0.45\textwidth}
\caption{Linear Zal+CN.}
\includegraphics[width=\textwidth, trim = 2.5cm 6cm 4cm 2cm, clip = true]
{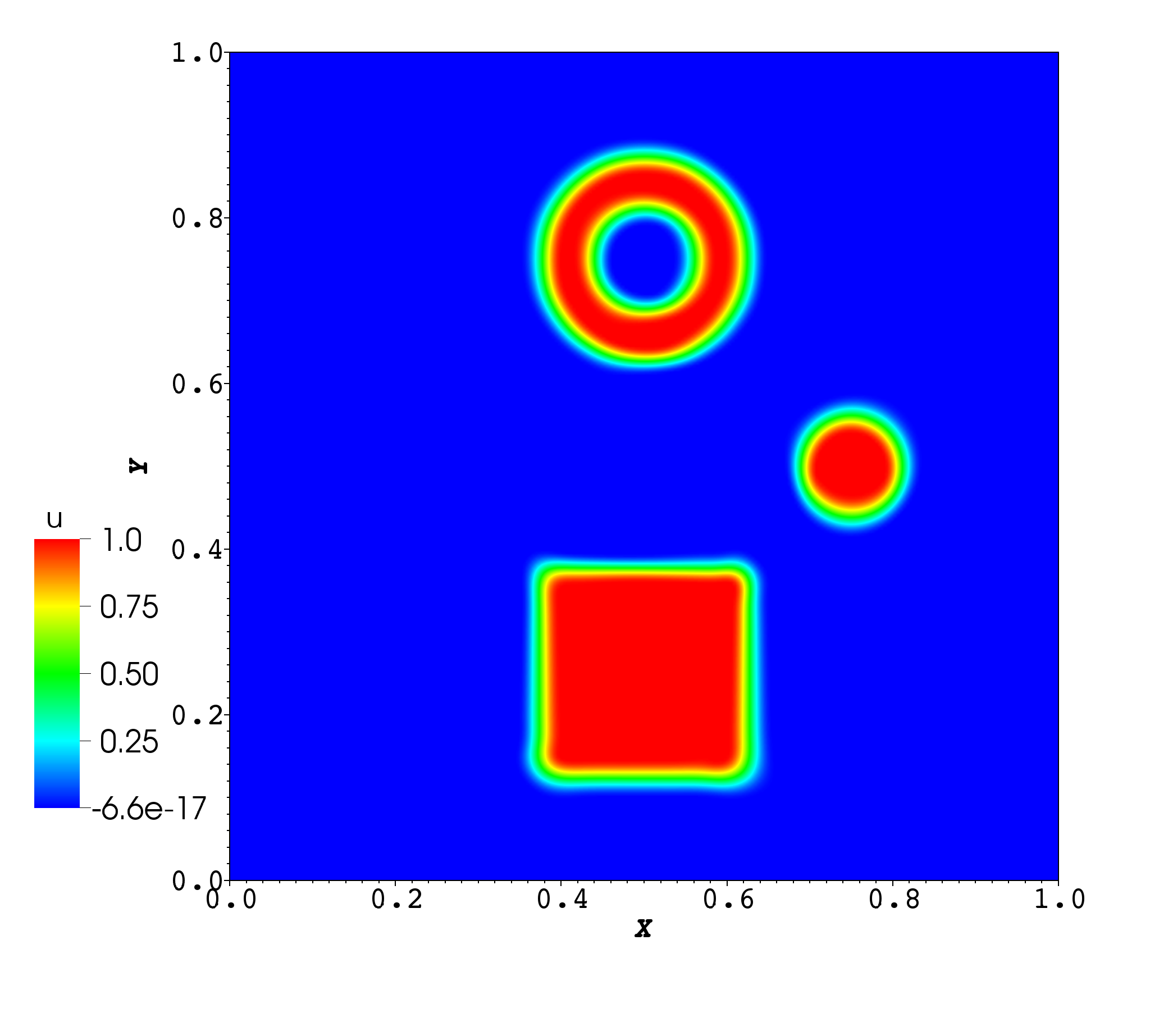}
\label{fig:linear_rotating_shapes_slice}
\end{subfigure}
\begin{subfigure}{0.45\textwidth}
\caption{Nonlinear Zal+CN.}
\includegraphics[width=\textwidth, trim = 2.5cm 6cm 4cm 2cm, clip = true]
{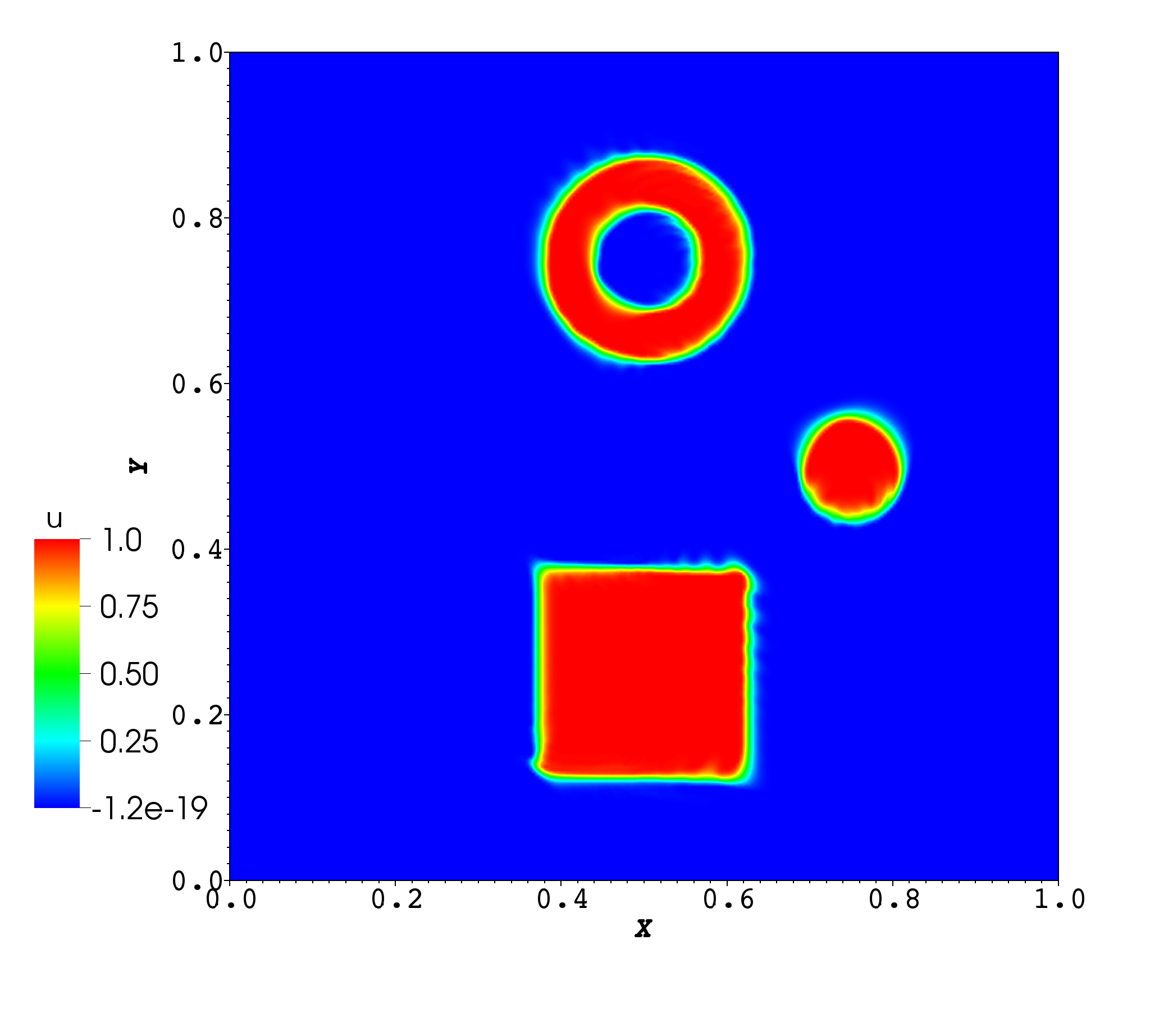}
\label{fig:nonlinear_rotating_shapes_slice}
\end{subfigure}
\begin{subfigure}{0.45\textwidth}
\caption{Zal+SSP.}
\includegraphics[width=\textwidth, trim = 2.5cm 6cm 4cm 2cm, clip = true]
{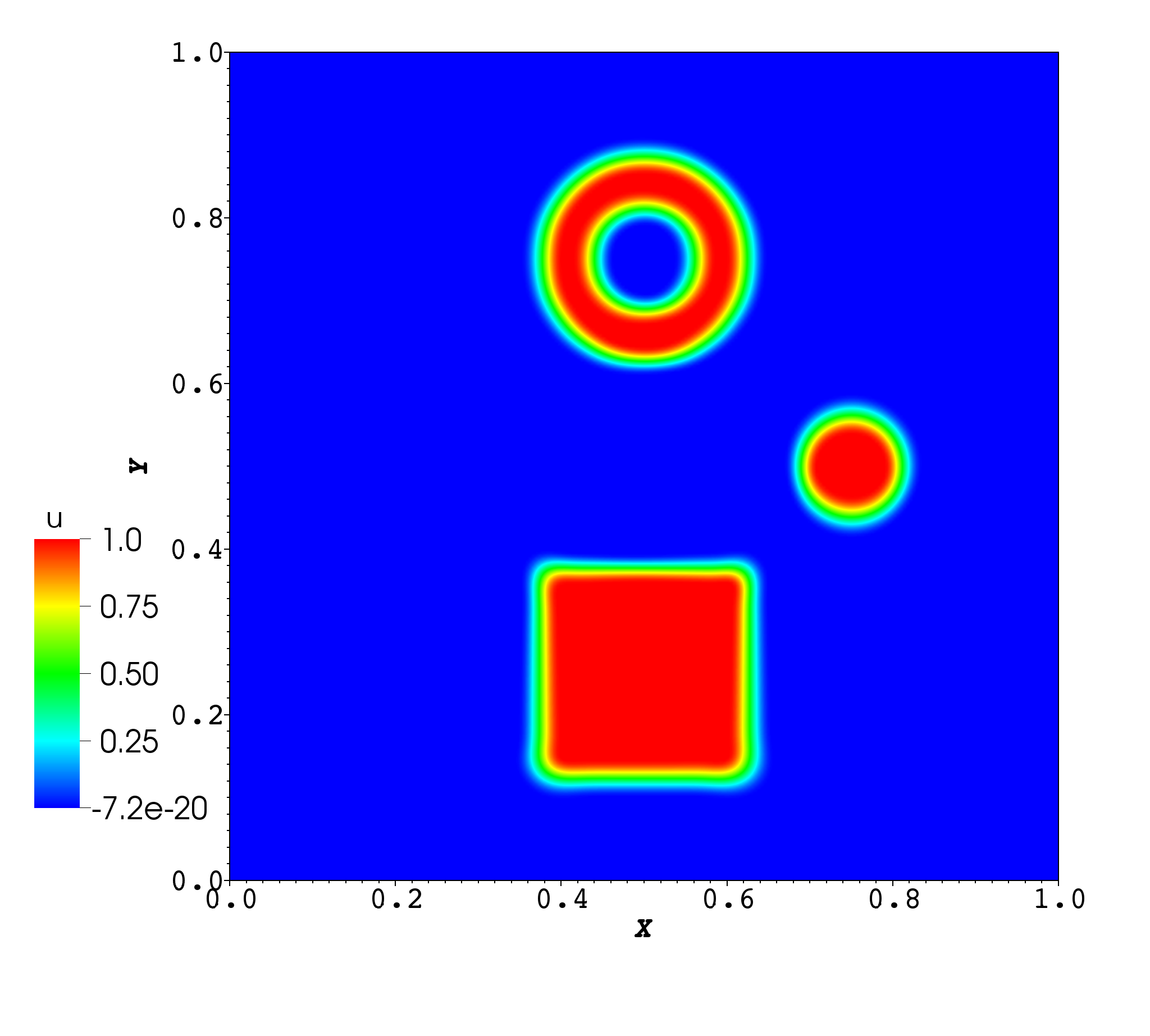}
\label{fig:zalesak_ssp}
\end{subfigure}
\begin{subfigure}{0.45\textwidth}
\caption{MC+CN.}
\includegraphics[width=\textwidth, trim = 2.5cm 6cm 4cm 2cm, clip = true]
{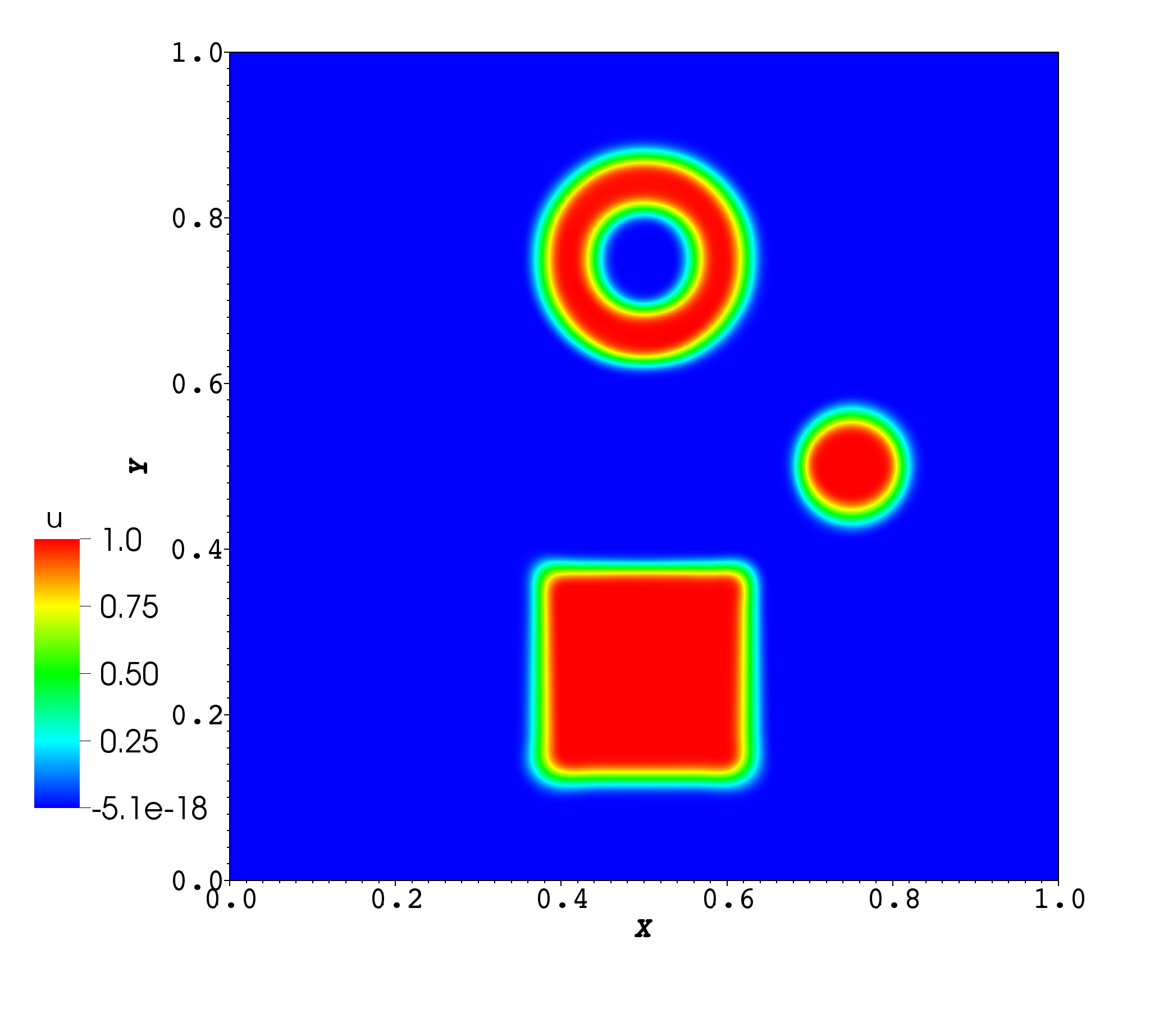}
\label{fig:monolithic_rotating_shapes_slice}
\end{subfigure}
\begin{subfigure}{0.45\textwidth}
\caption{MC+SSP.}
\includegraphics[width=\textwidth, trim = 2.5cm 6cm 4cm 2cm, clip = true]
{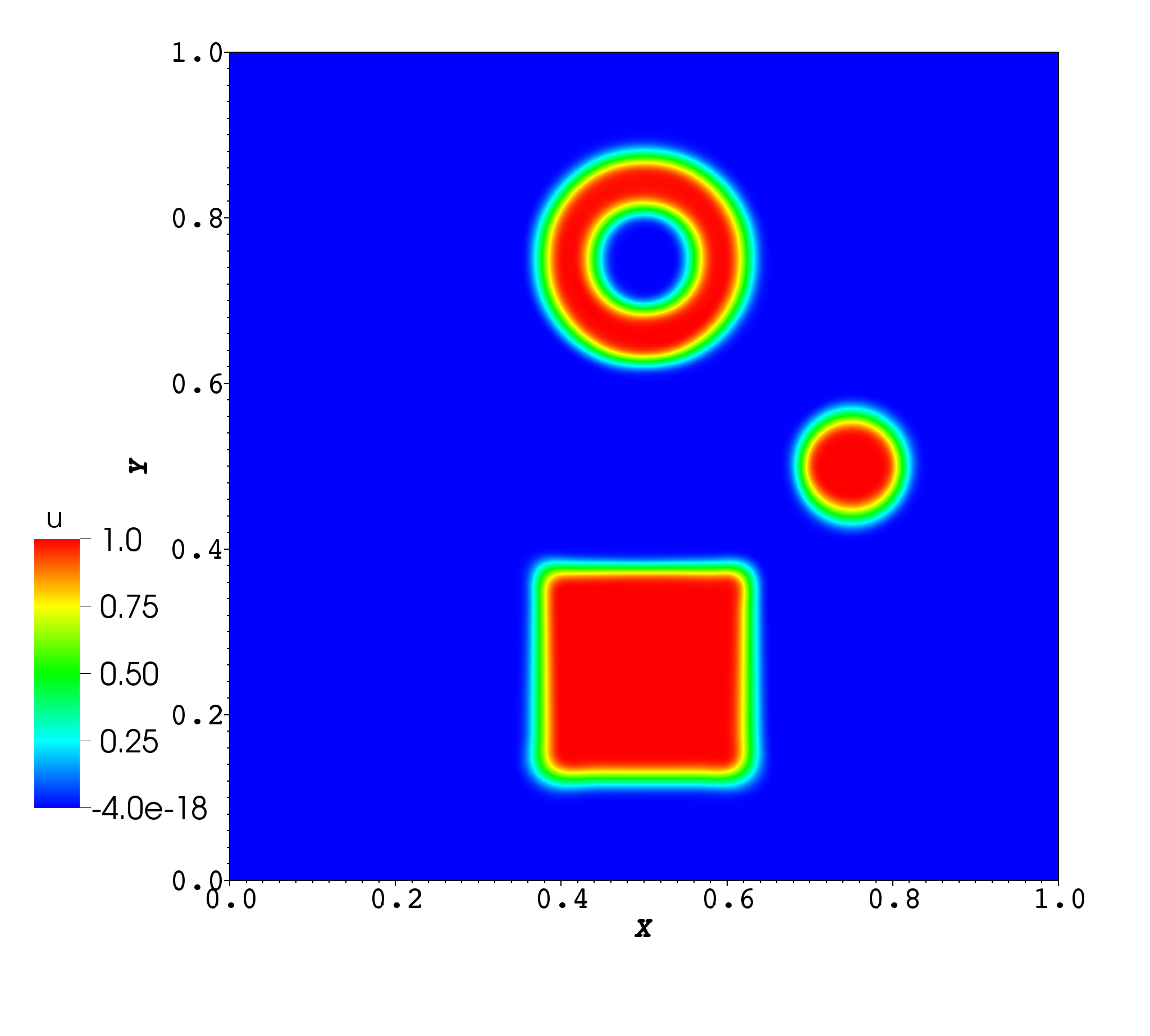}
\label{fig:monolithic_ssp}
\end{subfigure}
\label{fig:rotating_shapes_slice}
\end{figure}

The initial conditions and the solutions after a complete rotation in the plane $z = 0.5$ are shown in Figure~\ref{fig:rotating_shapes_slice}. The latter were calculated using all of the schemes mentioned in Section \ref{prelim:tcdr}.

Important criteria for assessing the quality of the solutions are the smearing of the layers and the size of the undershoots and overshoots.
This is compared in Figure~\ref{fig:rotating_shapes_slice} that shows the following:

The nonlinear Zal+CN smeared the solution the least, but the layers are uneven. This is probably due to our low limit on the number of the nonlinear iterations; the linear Zal+CN would also smear the layer very unevenly if we decreased $\mathtt{atol}$.

The other schemes smeared the solution quite uniformly;
but the schemes with the MC limiter or SSP produced slightly more uniform smearing than the linear Zal+CN.

None of the schemes produced undershoots or overshoots that were larger than the machine precision.

\begin{figure}[t!]
\centering
\caption{Example~\ref{ex:rotating_shapes}: Time evolution of the numerical solution and the true solution at the point $(0.5, 0.25, 0.5)^T$ (i.e., the initial center of the cube).}
\input{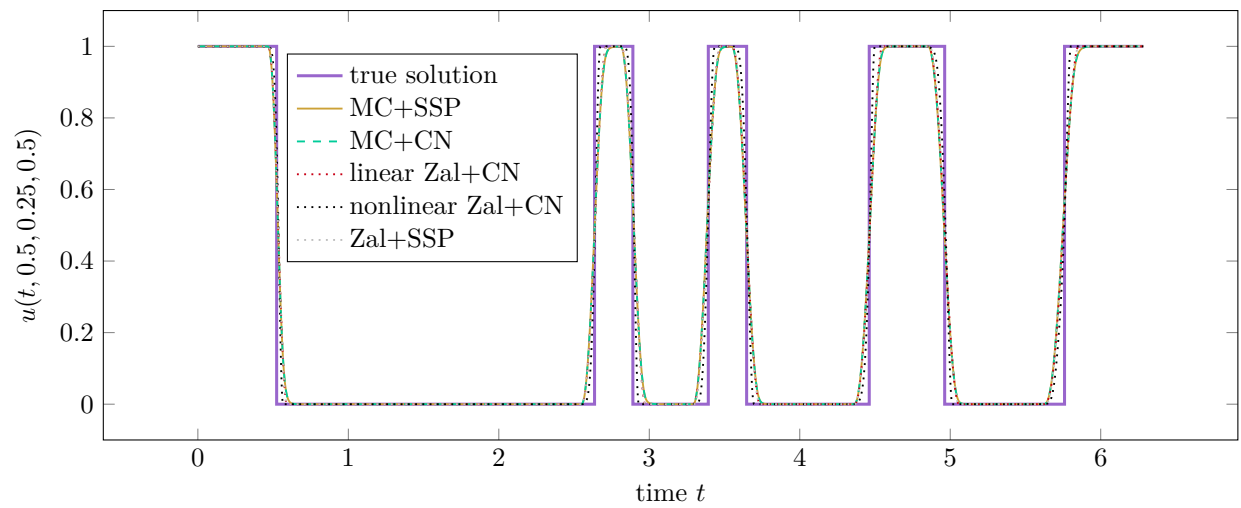}
\label{fig:rotating_shapes_points}
\end{figure}

Figure \ref{fig:rotating_shapes_points} compares the time evolution of the solution at the point $(0.5, 0.25, 0.5)^T$ (i.e., the initial center of the cube):
We can see that the schemes with the MC limiter produced essentially the same results. The results by the linear Zal+CN are almost identical to them, but they are a bit better in some parts: e.g., near the times 0.6 or 2.6.
The results of Zal+SSP are very similar to the linear Zal+CN.
The least diffusive results are produced by the nonlinear Zal+CN. However, all of the results are very similar.
The errors in $\left\Arrowvert \cdot \right\Arrowvert_2$ presented in Table~\ref{tab:rotating_shapes} confirm these observation.

\input{tables/rotating_shapes_L2_precTable.tex}

The largest difference between the schemes is in the computing times (Table~\ref{tab:rotating_shapes}):
The computing time by the nonlinear Zal+CN is by far the worst. The reason is that the rate of the decrease in the residual during the nonlinear loop steadily decreases in the course of the simulations.
(For MC+CN, on the contrary, this rate is roughly constant.) This property makes the scheme impractical for large problems and long simulations; and because of this, we decided to limit the number of the nonlinear iterations by 50.
Near the end of the simulation, the final norm of the residual in the nonlinear loop was around $10^{-13}$.

When comparing the nonlinear Zal+CN with the other schemes, we believe that the slightly higher precision of the solution does not compensate for the excessively long computing time at all.

One might be surprised that the explicit solvers are not the fastest ones because they do not require any (non)linear systems to be solved.
The reason is that one has to do the assembling twice in each time step (once for each SSP stage), and in our programs, this assembling is more expensive than solving the linear system.

Finally, Zal+SSP was faster than MC+SSP. This is because of the complexity of the assembling in the schemes with the MC limiter (see Remark \ref{rem:assembling_mc}, page \pageref{rem:assembling_mc}).

%% file: tables/rotating_shapes_L2_precTable.tex
\begin{table}
\centering
\caption{Example~\ref{ex:rotating_shapes}: Errors at final time measured in $\left\Arrowvert \cdot \right\Arrowvert_2$ and computing times.}
\begin{tabular}{c c c c c c}
 & \textbf{MC+SSP} & \textbf{MC+CN} & \textbf{Zal+CN} & \textbf{Zal+CN} & \textbf{Zal+SSP} \\  
 & & & \textbf{linear} & \textbf{nonlinear} & \\ \hline
\textbf{error} & 5.85e-2 & 5.85e-02 & 5.76e-2 & 4.58e-2 & 5.76e-2 \\ 
\textbf{computing time (hours)} & 29.9 & 106.6 & 17.7 & 254.4 & 22.7
\end{tabular}
\label{tab:rotating_shapes}
\end{table}

%% file: concentration_species2.tex
This example was proposed in \cite{JS08} and can be also found in \cite{JN12, AJ15}.
We use it to compare the efficiency of the solvers employed to solve the linear systems arising from the AFC schemes. We are also interested in the comparison between the Zalesak and the MC limiter.

\subsubsection*{Description of Problem}

This example models a typical situation where a species enters the domain at the entrance, travels through the domain while breeding and leaves it at the exit.

We consider~\eqref{eq:tcdr} for $\Omega = (0, 1)^3$ , $\varepsilon = 10^{-6}$, ${\bfb} = (1, -1/4, -1/8)^T$, $f = 0$ and
\begin{align}
c(x) = \begin{cases} 
        1 & \mbox{ if } \text{distance} (x,g) \leq 0.1,\\
        0 & \mbox{ otherwise},
      \end{cases}
\end{align}
where $g$ stands for the line segment with the endpoints $(0, 11/16, 11/16)^\text{T}$ and $(1, 7/16, 9/16)^\text{T}$.
These endpoints are the centers of the inlet $\Gamma_\text{in} = \lbrace 0 \rbrace \times [5/8,6/8] \times [5/8,6/8]$ and the outlet $\lbrace 1 \rbrace \times (3/8,4/8) \times (4/8,5/8) = \Gamma\nb$, where $\Gamma\db = \Gamma \setminus \Gamma\nb$.
Note the differences in the closeness of the above intervals that define the inlet and outlet.

We use the following boundary conditions:
At $\Gamma_\text{in}$, we prescribe
\begin{align}
u\db(t) = \begin{cases} 
        \sin (\pi t/2) & \mbox{ if } t\in [0,1],\\
        1 & \mbox{ if } t\in (1,2],\\
        \sin (\pi(t-1)/2) & \mbox{ if } t\in (2,3].
      \end{cases}
\end{align}
At $\Gamma\nb$, we prescribe $g\nb = 0$. At $\Gamma\db \setminus \Gamma_\text{in}$, we set $u\db=0$.

The initial condition is $u_0=0$, i.e., there are no species within the domain.

In the time interval $(0, 1)$, the inflow increases, and the injected species is transported to the outlet. Then, in the time interval $(1, 2)$, the inflow is constant, and the species reaches the outlet. Finally, in $(2, 3)$, the influx decreases.

The simulations were performed for the refinement levels 5--7, where the initial grid consisted of a single cube. The time step length was $\Delta t = 5 \cdot 10^{-3}$. Finally, we set $\mathtt{tol} = 10^{-14}$ in \eqref{eq:tol_formula} and $\mathtt{atol} = 10^{-16}$ for the \textsc{PETSc} solver.

Our choice of $\mathtt{tol}$ and $\mathtt{atol}$ means that the scheme MC+CN (unlike the other schemes) produces small undershoots, i.e., negative values of the magnitude larger than the machine precision.
This can be prevented by decreasing these tolerances. However, we did not decrease them because it would lead to excessive long computing times of the nonlinear Zal+CN (as in Section \ref{ex:rotating_shapes}). Instead, we replaced all negative components of the solution by zero in each time step.

We considered all of the solvers and preconditioners listed in Tables~\ref{tab:solvers_iterative}--\ref{tab:pc}. Simulations were performed with 4, 8, 16 and 32 processors for the refinement levels 5--7.
However, in the sequential case, we performed the simulations only for the refinement level~5. Note that we do not list results of the parallel solvers for BCGS + Jac below because the \textsc{PETSc} solver always crashed for this combination.

\subsubsection*{Discussion of Results}

Tables~\ref{tab:nonlinear_mpi_5}--\ref{tab:nonlinear_mpi_7} (pages \pageref{tab:nonlinear_mpi_5}--\pageref{tab:nonlinear_mpi_7}),  \ref{tab:linear_mpi_5}--\ref{tab:linear_mpi_7} (pages \pageref{tab:linear_mpi_5}--\pageref{tab:linear_mpi_7}) and \ref{tab:MC_mpi_5}--\ref{tab:MC_mpi_7} (pages \pageref{tab:MC_mpi_5}--\pageref{tab:MC_mpi_7}) show the dependence of the computing time on the number of processors for the nonlinear and linear Zal+CN and for MC+CN, respectively. We believe that the reduction in computation time by doubling the number of processors is acceptable.

\input{tables/concentration_species/nonlinear_5_mpi2_table.tex}
\input{tables/concentration_species/nonlinear_6_mpi2_table.tex}
\input{tables/concentration_species/nonlinear_7_mpi2_table.tex}

\input{tables/concentration_species/linear_5_mpi2_table.tex}
\input{tables/concentration_species/linear_6_mpi2_table.tex}
\input{tables/concentration_species/linear_7_mpi2_table.tex}

\input{tables/concentration_species/monolithic_5_table.tex}
\input{tables/concentration_species/monolithic_6_table.tex}
\input{tables/concentration_species/monolithic_7_table.tex}

We can clearly see that, generally, the best solver is FGMRES, and the worst solver is BCGS. The order of preconditioners is (from the best one): Jac, SOR, BJac, ASM and MG, where ASM performs much worse than the first three preconditioners, and MG performs much worse than ASM.

As expected, the computing times for the nonlinear Zal+CN are several times longer than those for MC+CN.

In the above discussion, we compared only the parallel algorithms because we expected them to perform better than their sequential counterparts for the considered sizes of the problems.
The correctness of this conjecture is illustrated by Tables \ref{tab;nonlinear_sequential}--\ref{tab;monolithic_sequential} (pages \pageref{tab;nonlinear_sequential}--\pageref{tab;monolithic_sequential}) in which we listed the computing times for sequential and parallel (4 processors) solvers on the coarse refinement level. The parallel versions are clearly better.

These tables also show the computing time for the direct solvers LU and UMFPACK. Although UMFPACK is clearly competitive in the sequential case, the computing time for parallel versions are by far worse than those for the iterative solvers.
Therefore, we decided not to include them in our further comparisons. Also note the interesting fact that the computing time with UMFPACK actually shots up when running the solver in parallel.

\input{tables/concentration_species/nonlinear_sequential}
\input{tables/concentration_species/linear_sequential}
\input{tables/concentration_species/monolithic_sequential}

\begin{figure}[t!]
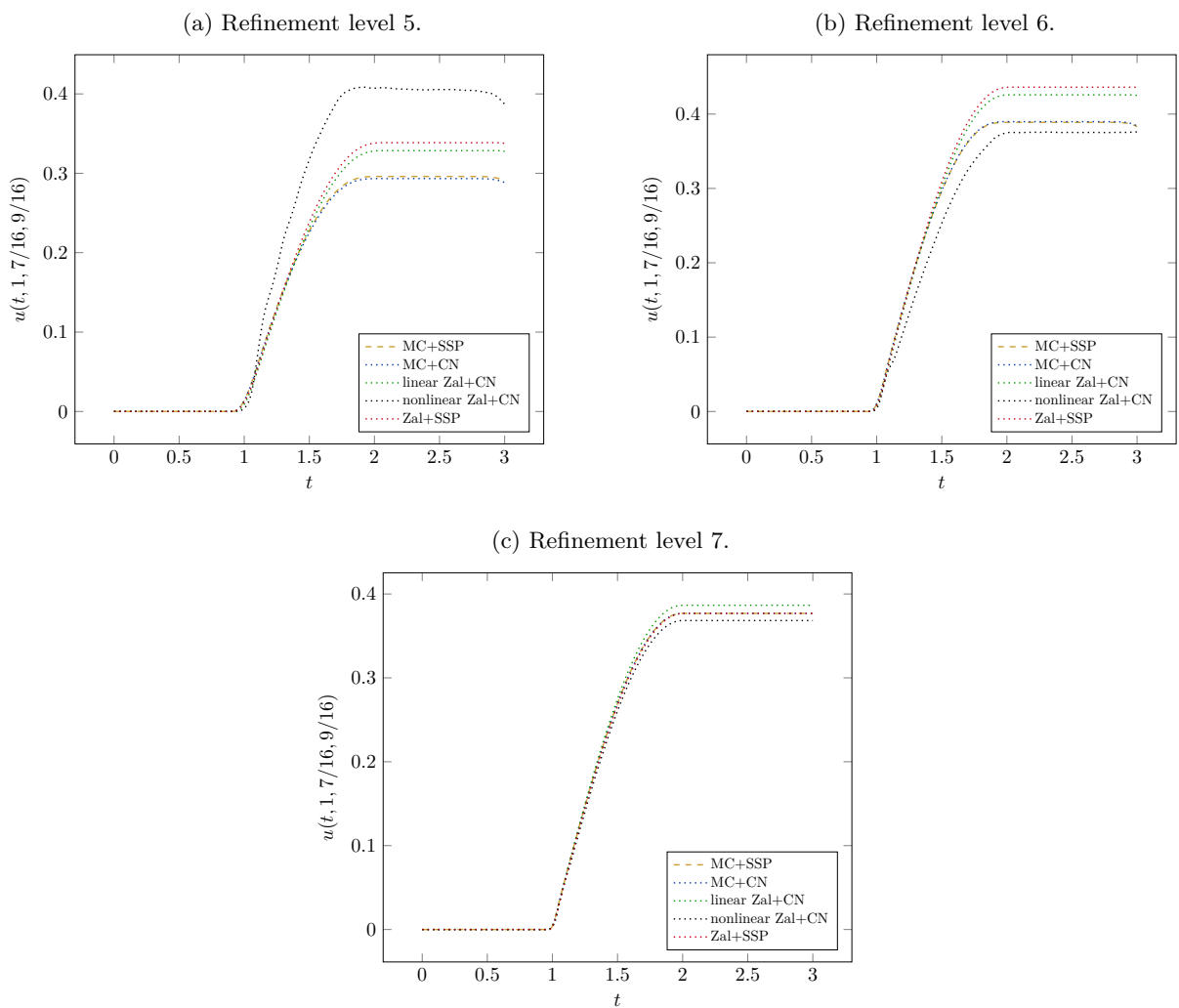

\centering
\caption{Example~\ref{ex:concentration_species2}: Time evolution of $u$ at $(1,7/16,9/16)^\text{T}$. The MC limiter produces essentially the same result for both time discretizations.}
\begin{subfigure}{0.49\textwidth}
\caption{Refinement level 5.}
\input{graphs/concentration_species/concentration_species_point_level_5b.tex}
\label{fig:concentration_point_level_5}
\end{subfigure}\hspace*{1em}
\begin{subfigure}{0.49\textwidth}
\caption{Refinement level 6.}
\input{graphs/concentration_species/concentration_species_point_level_6b.tex}
\label{fig:concentration_point_level_6}
\end{subfigure}
\par\bigskip
\begin{subfigure}{0.49\textwidth}
\caption{Refinement level 7.}
\input{graphs/concentration_species/concentration_species_point_level_7b.tex}
\label{fig:concentration_point_level_7}
\end{subfigure}
\label{fig:concentration_species_point}
\end{figure}

As a measure of accuracy, the authors of \cite{JN12} proposed to compare the time evolution of the solution at the center of the outlet, i.e., at the point $(1, 7/16, 9/16)^\text{T}$. This is done in Figure~\ref{fig:concentration_species_point}.

The implicit schemes are represented only by the combination FGMRES + SOR; the other combinations produced very similar results.
Contrary to the implicit schemes, the explicit schemes were used with the time step $\Delta t = 10^{-3}$ (because of the CFL condition).
The curves representing the time evolutions probably converge to some curve, but it is not clear which one is better. Note that the schemes with the MC limiter produced essentially the same results.

%% file: tables/concentration_species/nonlinear_5_mpi2_table.tex
\begin{table}[t!]
\centering
\caption{Example~\ref{ex:concentration_species2}: Computing time in seconds of the solver with the nonlinear Zal+CN for refinement level~5.}
\begin{tabular}{c  c  c c c c}
\multicolumn{2}{c}{ }  & \multicolumn{4}{c }{\textbf{time (sec)}} \\
\textbf{Solver} & \textbf{PC} & NP = 4 & NP = 8 & NP = 16 & NP = 32 \\ \hline
\multirow{4}{*}{BCGS} & BJac & 2079 & 1222 & 760 & 512 \\ \cline{2-6}
 & SOR & 1986 & 1167 & 726 & 496 \\ \cline{2-6}
 & ASM & 2398 & 1452 & 987 & 680 \\ \cline{2-6}
 & MG & 2540 & 1476 & 1081 & 895 \\ \hline
\multirow{5}{*}{LGMRES} & BJac & 2079 & 1229 & 763 & 518 \\ \cline{2-6}
 & SOR & 1989 & 1178 & 726 & 499 \\ \cline{2-6}
 & ASM & 2378 & 1444 & 986 & 678 \\ \cline{2-6}
 & Jac & 1901 & 1130 & 698 & 485 \\ \cline{2-6}
 & MG & 2534 & 1512 & 1115 & 955 \\ \hline
\multirow{5}{*}{FGMRES} & BJac & 2003 & 1185 & 737 & 499 \\ \cline{2-6}
 & SOR & 1907 & 1126 & 706 & 476 \\ \cline{2-6}
 & ASM & 2290 & 1398 & 950 & 650 \\ \cline{2-6}
 & Jac & 1860 & 1107 & 684 & 471 \\ \cline{2-6}
 & MG & 2270 & 1402 & 1067 & 957 \\
\end{tabular}
\label{tab:nonlinear_mpi_5}
\end{table}

%% file: tables/concentration_species/nonlinear_6_mpi2_table.tex
\begin{table}[t!]
\centering
\caption{Example~\ref{ex:concentration_species2}: Computing time in seconds of the solver with the nonlinear Zal+CN for refinement level~6.}
\begin{tabular}{c  c  c c c c}
\multicolumn{2}{c}{ } & \multicolumn{4}{c}{\textbf{time (sec)}} \\ 
\textbf{Solver} & \textbf{PC} & NP = 4 & NP = 8 & NP = 16 & NP = 32 \\ \hline
\multirow{4}{*}{BCGS} & BJac & 12488 & 6670 & 3772 & 2248 \\ \cline{2-6}
 & SOR & 12004 & 6314 & 3596 & 2124 \\ \cline{2-6}
 & ASM & 14144 & 7555 & 4463 & 2760 \\ \cline{2-6}
 & MG & 17302 & 8763 & 5187 & 3022 \\ \hline
\multirow{5}{*}{LGMRES} & BJac & 12459 & 6626 & 3787 & 2254 \\ \cline{2-6}
 & SOR & 11948 & 6389 & 3612 & 2128 \\ \cline{2-6}
 & ASM & 14077 & 7528 & 4448 & 2736 \\ \cline{2-6}
 & Jac & 11427 & 6115 & 3445 & 2046 \\ \cline{2-6}
 & MG & 17225 & 8786 & 5239 & 3081 \\ \hline
\multirow{5}{*}{FGMRES} & BJac & 11682 & 6232 & 3515 & 2090 \\ \cline{2-6}
 & SOR & 11068 & 5939 & 3354 & 2002 \\ \cline{2-6}
 & ASM & 13297 & 7162 & 4217 & 2561 \\ \cline{2-6}
 & Jac & 10846 & 5873 & 3280 & 1964 \\ \cline{2-6}
 & MG & 14636 & 7603 & 4470 & 2689 \\
\end{tabular}
\label{tab:nonlinear_mpi_6}
\end{table}

%% file: tables/concentration_species/nonlinear_7_mpi2_table.tex
\begin{table}[t!]
\centering
\caption{Example~\ref{ex:concentration_species2}: Computing time in seconds of the solver with the nonlinear Zal+CN for refinement level~7.}
\begin{tabular}{c c c c c c}
\multicolumn{2}{c}{ } & \multicolumn{4}{c}{\textbf{time (sec)}} \\ 
\textbf{Solver} & \textbf{PC} & NP = 4 & NP = 8 & NP = 16 & NP = 32 \\\hline
\multirow{4}{*}{BCGS} & BJac & 52053 & 29114 & 15469 & 8941 \\ \cline{2-6}
 & SOR & 51732 & 29042 & 15504 & 9209 \\ \cline{2-6}
 & ASM & 58398 & 32978 & 17694 & 10356 \\ \cline{2-6}
 & MG & 89549 & 54378 & 31340 & 16829 \\ \hline
\multirow{5}{*}{LGMRES} & BJac & 52176 & 29198 & 15598 & 9023 \\ \cline{2-6}
 & SOR & 52203 & 29426 & 15508 & 9293 \\ \cline{2-6}
 & ASM & 58760 & 32953 & 17689 & 10361 \\ \cline{2-6}
 & Jac & 50720 & 27788 & 14815 & 9035 \\ \cline{2-6}
 & MG & 91239 & 47312 & 26151 & 17150 \\ \hline
\multirow{5}{*}{FGMRES} & BJac & 47523 & 26532 & 14117 & 8141 \\ \cline{2-6}
 & SOR & 46234 & 25695 & 13666 & 8010 \\ \cline{2-6}
 & ASM & 54173 & 30610 & 16409 & 9519 \\ \cline{2-6}
 & Jac & 45626 & 24908 & 13356 & 7846 \\ \cline{2-6}
 & MG & 64826 & 36046 & 20632 & 12251 \\
\end{tabular}
\label{tab:nonlinear_mpi_7}
\end{table}

%% file: tables/concentration_species/linear_5_mpi2_table.tex
\begin{table}[t!]
\centering
\caption{Example~\ref{ex:concentration_species2}: Computing time in seconds of the solver with the linear Zal+CN for refinement level~5.}
\begin{tabular}{c c c c c c}
\multicolumn{2}{c}{ } & \multicolumn{4}{c}{\textbf{time (sec)}} \\ 
\textbf{Solver} & \textbf{PC} & NP = 4 & NP = 8 & NP = 16 & NP = 32 \\ \hline
\multirow{4}{*}{BCGS} & BJac & 82 & 46 & 30 & 21 \\ \cline{2-6}
 & SOR & 82 & 46 & 30 & 20 \\ \cline{2-6}
 & ASM & 87 & 48 & 33 & 21 \\ \cline{2-6}
 & MG & 87 & 48 & 33 & 22 \\ \hline
\multirow{5}{*}{LGMRES} & BJac & 82 & 46 & 31 & 20 \\ \cline{2-6}
 & SOR & 82 & 45 & 30 & 20 \\ \cline{2-6}
 & ASM & 86 & 48 & 33 & 21 \\ \cline{2-6}
 & Jac & 81 & 44 & 30 & 19 \\ \cline{2-6}
 & MG & 89 & 49 & 32 & 22 \\ \hline
\multirow{5}{*}{FGMRES} & BJac & 81 & 45 & 30 & 19 \\ \cline{2-6}
 & SOR & 82 & 45 & 30 & 20 \\ \cline{2-6}
 & ASM & 86 & 48 & 32 & 21 \\ \cline{2-6}
 & Jac & 81 & 45 & 30 & 21 \\ \cline{2-6}
 & MG & 85 & 48 & 32 & 21 \\ 
\end{tabular}
\label{tab:linear_mpi_5}
\end{table}

%% file: tables/concentration_species/linear_6_mpi2_table.tex
\begin{table}[t!]
\centering
\caption{Example~\ref{ex:concentration_species2}: Computing time in seconds of the solver with the linear Zal+CN for refinement level~6.}
\begin{tabular}{c c c c c c}
\multicolumn{2}{c}{ }  & \multicolumn{4}{c}{\textbf{time (sec)}} \\
\textbf{Solver} & \textbf{PC} & NP = 4 & NP = 8 & NP = 16 & NP = 32 \\ \hline
\multirow{4}{*}{BCGS} & BJac & 615 & 339 & 194 & 114 \\ \cline{2-6}
 & SOR & 616 & 331 & 194 & 115 \\ \cline{2-6}
 & ASM & 626 & 343 & 202 & 121 \\ \cline{2-6}
 & MG & 665 & 362 & 215 & 125 \\ \hline
\multirow{5}{*}{LGMRES} & BJac & 624 & 342 & 198 & 116 \\ \cline{2-6}
 & SOR & 617 & 336 & 194 & 115 \\ \cline{2-6}
 & ASM & 650 & 351 & 207 & 123 \\ \cline{2-6}
 & Jac & 615 & 328 & 189 & 112 \\ \cline{2-6}
 & MG & 701 & 370 & 214 & 128 \\ \hline
\multirow{5}{*}{FGMRES} & BJac & 598 & 326 & 192 & 114 \\ \cline{2-6}
 & SOR & 604 & 329 & 188 & 113 \\ \cline{2-6}
 & ASM & 633 & 338 & 198 & 120 \\ \cline{2-6}
 & Jac & 605 & 327 & 188 & 112 \\ \cline{2-6}
 & MG & 650 & 351 & 202 & 122 \\ 
\end{tabular}
\label{tab:linear_mpi_6}
\end{table}

%% file: tables/concentration_species/linear_7_mpi2_table.tex
\begin{table}[t!]
\centering
\caption{Example~\ref{ex:concentration_species2}: Computing time in seconds of the solver with the linear Zal+CN for refinement level~7.}
\begin{tabular}{c c c c c c}
\multicolumn{2}{c}{ }  & \multicolumn{4}{c}{\textbf{time (sec)}} \\ 
\textbf{Solver} & \textbf{PC} & NP = 4 & NP = 8 & NP = 16 & NP = 32 \\\hline
\multirow{4}{*}{BCGS} & BJac & 4905 & 2612 & 1406 & 794 \\ \cline{2-6}
 & SOR & 4870 & 2600 & 1433 & 808 \\ \cline{2-6}
 & ASM & 4974 & 2696 & 1468 & 822 \\ \cline{2-6}
 & MG & 5747 & 3142 & 1797 & 1019 \\ \hline
\multirow{5}{*}{LGMRES} & BJac & 5045 & 2664 & 1440 & 813 \\ \cline{2-6}
 & SOR & 5097 & 2677 & 1467 & 846 \\ \cline{2-6}
 & ASM & 5198 & 2820 & 1554 & 878 \\ \cline{2-6}
 & Jac & 4827 & 2597 & 1381 & 803 \\ \cline{2-6}
 & MG & 5971 & 3224 & 1882 & 1020 \\ \hline
\multirow{5}{*}{FGMRES} & BJac & 4831 & 2606 & 1417 & 796 \\ \cline{2-6}
 & SOR & 4823 & 2566 & 1387 & 794 \\ \cline{2-6}
 & ASM & 4976 & 2696 & 1481 & 834 \\ \cline{2-6}
 & Jac & 4706 & 2565 & 1397 & 783 \\ \cline{2-6}
 & MG & 5465 & 2981 & 1650 & 935 \\
\end{tabular}
\label{tab:linear_mpi_7}
\end{table}

%% file: tables/concentration_species/monolithic_5_table.tex
\begin{table}[t!]
\centering
\caption{Example~\ref{ex:concentration_species2}: Computing time in seconds of the solver with MC+CN for refinement level~5.}
\begin{tabular}{c c c c c c}
\multicolumn{2}{c}{ } & \multicolumn{4}{c}{\textbf{time (sec)}} \\ 
\textbf{Solver} & \textbf{PC} & NP = 4 & NP = 8 & NP = 16 & NP = 32 \\\hline
\multirow{4}{*}{BCGS} & BJac & 184 & 102 & 66 & 42 \\ \cline{2-6}
 & SOR & 177 & 99 & 64 & 42 \\ \cline{2-6}
 & ASM & 196 & 110 & 75 & 49 \\ \cline{2-6}
 & MG & 201 & 113 & 73 & 48 \\ \hline
\multirow{5}{*}{LGMRES} & BJac & 184 & 102 & 66 & 43 \\ \cline{2-6}
 & SOR & 178 & 99 & 64 & 42 \\ \cline{2-6}
 & ASM & 195 & 110 & 76 & 50 \\ \cline{2-6}
 & Jac & 173 & 97 & 63 & 41 \\ \cline{2-6}
 & MG & 201 & 111 & 74 & 50 \\ \hline
\multirow{5}{*}{FGMRES} & BJac & 177 & 99 & 65 & 42 \\ \cline{2-6}
 & SOR & 175 & 98 & 64 & 41 \\ \cline{2-6}
 & ASM & 198 & 109 & 74 & 48 \\ \cline{2-6}
 & Jac & 173 & 96 & 63 & 41 \\ \cline{2-6}
 & MG & 212 & 107 & 70 & 46 \\
\end{tabular}
\label{tab:MC_mpi_5}
\end{table}

%% file: tables/concentration_species/monolithic_6_table.tex
\begin{table}[t!]
\centering
\caption{Example~\ref{ex:concentration_species2}: Computing time in seconds of the solver with MC+CN for refinement level~6.}
\begin{tabular}{c c c c c c}
\multicolumn{2}{c}{ } & \multicolumn{4}{c}{\textbf{time (sec)}} \\ 
\textbf{Solver} & \textbf{PC} & NP = 4 & NP = 8 & NP = 16 & NP = 32 \\\hline
\multirow{4}{*}{BCGS} & BJac & 1349 & 736 & 416 & 247 \\ \cline{2-6}
 & SOR & 1329 & 717 & 407 & 240 \\ \cline{2-6}
 & ASM & 1429 & 777 & 453 & 276 \\ \cline{2-6}
 & MG & 1604 & 852 & 490 & 284 \\ \hline
\multirow{5}{*}{LGMRES} & BJac & 1364 & 731 & 420 & 249 \\ \cline{2-6}
 & SOR & 1332 & 721 & 411 & 243 \\ \cline{2-6}
 & ASM & 1448 & 790 & 460 & 276 \\ \cline{2-6}
 & Jac & 1310 & 732 & 403 & 239 \\ \cline{2-6}
 & MG & 1608 & 857 & 496 & 287 \\ \hline
\multirow{5}{*}{FGMRES} & BJac & 1320 & 706 & 404 & 239 \\ \cline{2-6}
 & SOR & 1281 & 689 & 396 & 234 \\ \cline{2-6}
 & ASM & 1396 & 762 & 444 & 267 \\ \cline{2-6}
 & Jac & 1271 & 694 & 393 & 234 \\ \cline{2-6}
 & MG & 1447 & 779 & 450 & 264 \\
\end{tabular}
\label{tab:MC_mpi_6}
\end{table}

%% file: tables/concentration_species/monolithic_7_table.tex
\begin{table}[t!]
\centering
\caption{Example~\ref{ex:concentration_species2}: Computing time in seconds of the solver with MC+CN for refinement level~7.}
\begin{tabular}{c c c c c c}
\multicolumn{2}{c}{ } & \multicolumn{4}{c}{\textbf{time (sec)}} \\ 
\textbf{Solver} & \textbf{PC} & NP = 4 & NP = 8 & NP = 16 & NP = 32 \\\hline
\multirow{4}{*}{BCGS} & BJac & 10127 & 5449 & 2938 & 1677 \\ \cline{2-6}
 & SOR & 10174 & 5428 & 2934 & 1702 \\ \cline{2-6}
 & ASM & 10842 & 5774 & 3141 & 1778 \\ \cline{2-6}
 & MG & 13790 & 7209 & 3898 & 2360 \\ \hline
\multirow{5}{*}{LGMRES} & BJac & 10374 & 5566 & 2985 & 1702 \\ \cline{2-6}
 & SOR & 10446 & 5569 & 2996 & 1739 \\ \cline{2-6}
 & ASM & 11072 & 5926 & 3220 & 1845 \\ \cline{2-6}
 & Jac & 10311 & 5387 & 2924 & 1700 \\ \cline{2-6}
 & MG & 14058 & 7295 & 3954 & 2387 \\ \hline
\multirow{5}{*}{FGMRES} & BJac & 9924 & 5278 & 2837 & 1602 \\ \cline{2-6}
 & SOR & 9762 & 5490 & 2812 & 1603 \\ \cline{2-6}
 & ASM & 10554 & 5661 & 3069 & 1742 \\ \cline{2-6}
 & Jac & 9692 & 5144 & 2779 & 1593 \\ \cline{2-6}
 & MG & 11987 & 6267 & 3386 & 2024 \\
\end{tabular}
\label{tab:MC_mpi_7}
\end{table}

%% file: tables/concentration_species/nonlinear_sequential.tex
\begin{table}[t!]
\centering
\caption{Example~\ref{ex:concentration_species2}: Nonlinear Zal+CN: Computing time in seconds for $\text{NP} = 1, 4$ and refinement level 5.}
\begin{tabular}{ c c  c c c c c c}
\multicolumn{2}{c}{} & \multicolumn{2}{c}{\textbf{time (sec)}}&\multicolumn{2}{c}{} &\multicolumn{2}{c}{\textbf{time (sec)}}\\ 
\textbf{Solver} & \textbf{PC} & NP = 1 & NP = 4 & \textbf{Solver} & \textbf{PC} & NP = 1 & NP = 4\\ \hline
\multirow{5}{*}{BCGS} & BJac & 5740 & 2079 & \multirow{5}{*}{FGMRES} & BJac & 5288 & 2003 \\ \cline{2-4} \cline{6-8}
 & SOR & 5355 & 1986 &  & SOR & 4958 & 1907 \\ \cline{2-4} \cline{6-8}
 & ASM & 6271 & 2398 &  & ASM & 5821 & 2290 \\ \cline{2-4} \cline{6-8}
 & Jac & 5224 & --- &  & Jac & 4958 & 1860 \\ \cline{2-4} \cline{6-8}
 & MG & 6618 & 2540 &  & MG & 6013 & 2270 \\ \hline
\multirow{5}{*}{LGMRES} & BJac & 5915 & 2079 & LU & --- & 85567 & 48247 \\ \cline{2-4} \cline{5-8}
 & SOR & 5248 & 1989 & UMFPACK & --- & 5073 & 48987 \\ \cline{2-4} \cline{5-8}
 & ASM & 6325 & 2378 & & & & \\ \cline{2-4} \cline{5-8}
 & Jac & 5077 & 1901 & & & & \\ \cline{2-4} \cline{5-8}
 & MG & 6599 & 2534 & & & & \\
\end{tabular}
\label{tab;nonlinear_sequential}
\end{table}

%% file: tables/concentration_species/linear_sequential.tex
\begin{table}[t!]
\centering
\caption{Example~\ref{ex:concentration_species2}: Linear Zal+CN: Computing time in seconds for $\text{NP} = 1, 4$ and refinement level 5.}
\begin{tabular}{ c c c c c c c c }
\multicolumn{2}{ c}{} & \multicolumn{2}{c}{\textbf{time (sec)}}&\multicolumn{2}{c}{} &\multicolumn{2}{c}{\textbf{time (sec)}}\\
\textbf{Solver} & \textbf{PC} & NP = 1 & NP = 4 & \textbf{Solver} & \textbf{PC} & NP = 1 & NP = 4\\ \hline
\multirow{5}{*}{BCGS} & BJac & 252 & 82  & \multirow{5}{*}{FGMRES} & BJac  & 249 & 81 \\ \cline{2-4} \cline{6-8}
 & SOR & 248 & 82  &  & SOR   & 247 & 82 \\ \cline{2-4} \cline{6-8}
 & ASM & 253 & 87  &  & ASM   & 254 & 86 \\ \cline{2-4} \cline{6-8}
 & Jac & 248  & ---  &  & Jac   & 248 & 81 \\ \cline{2-4} \cline{6-8}
 & MG  & 261 & 87  &  & MG    & 255 & 85 \\ \hline
\multirow{5}{*}{LGMRES} & BJac  & 252 & 82  & LU & --- & 1009 & 513   \\ \cline{2-4} \cline{5-8}
 & SOR   & 249 & 82  & UMFPACK & --- & 248 & 516   \\ \cline{2-4} \cline{5-8}
 & ASM   & 255 & 86 & & & &\\ \cline{2-4} \cline{5-8}
 & Jac   & 250 & 81 & & & &\\ \cline{2-4} \cline{5-8}
 & MG    & 261 & 89 & & & &\\
\end{tabular}
\label{tab;linear_sequential}
\end{table}

%% file: tables/concentration_species/monolithic_sequential.tex
\begin{table}[t!]
\centering
\caption{Example~\ref{ex:concentration_species2}: MC+CN: Computing time in seconds for $\text{NP} = 1, 4$ and refinement level 5.}
\begin{tabular}{c c c c c c c c }
\multicolumn{2}{c}{} & \multicolumn{2}{c}{\textbf{time (sec)}}&\multicolumn{2}{c}{} &\multicolumn{2}{c}{\textbf{time (sec)}}\\ 
\textbf{Solver} & \textbf{PC} & NP = 1 & NP = 4 & \textbf{Solver} & \textbf{PC} & NP = 1 & NP = 4\\ \hline
\multirow{5}{*}{BCGS} & BJac & 553 & 184  & \multirow{5}{*}{FGMRES} & BJac  & 543 & 177 \\ \cline{2-4} \cline{6-8}
 & SOR & 543 & 177  &  & SOR   & 532 & 175 \\ \cline{2-4} \cline{6-8}
 & ASM & 585 & 196  &  & ASM   & 560 & 198 \\ \cline{2-4} \cline{6-8}
 & Jac & 538  & ---  &  & Jac  & 534 & 173 \\ \cline{2-4} \cline{6-8}
 & MG  & 593 & 201  &  & MG    & 564 & 212 \\ \hline
\multirow{5}{*}{LGMRES} & BJac & 558 & 184  & LU & --- & 3750 & 2066   \\ \cline{2-4} \cline{5-8}
 & SOR   & 541 & 178  & UMFPACK & --- & 513 & 2078   \\ \cline{2-4} \cline{5-8}
 & ASM   & 571 & 195 & & & &\\ \cline{2-4} \cline{5-8}
 & Jac   & 538 & 173 & & & &\\ \cline{2-4} \cline{5-8}
 & MG    & 590 & 201 & & & &\\
\end{tabular}
\label{tab;monolithic_sequential}
\end{table}

%% file: example_non_constant_convection2.tex
This example was proposed in \cite{BJKR18}. We use it to compare the computing times of different solvers.

\subsubsection*{Description of Problem}

\begin{figure}[t!]
\caption{Example~\ref{ex:non_constant_convection2}: Isosurface for $u=0.05$ of the solution computed using the MU limiter for the level~5, and sketch of the coarsest grid (level~0).}
\includegraphics[width=0.5\textwidth]{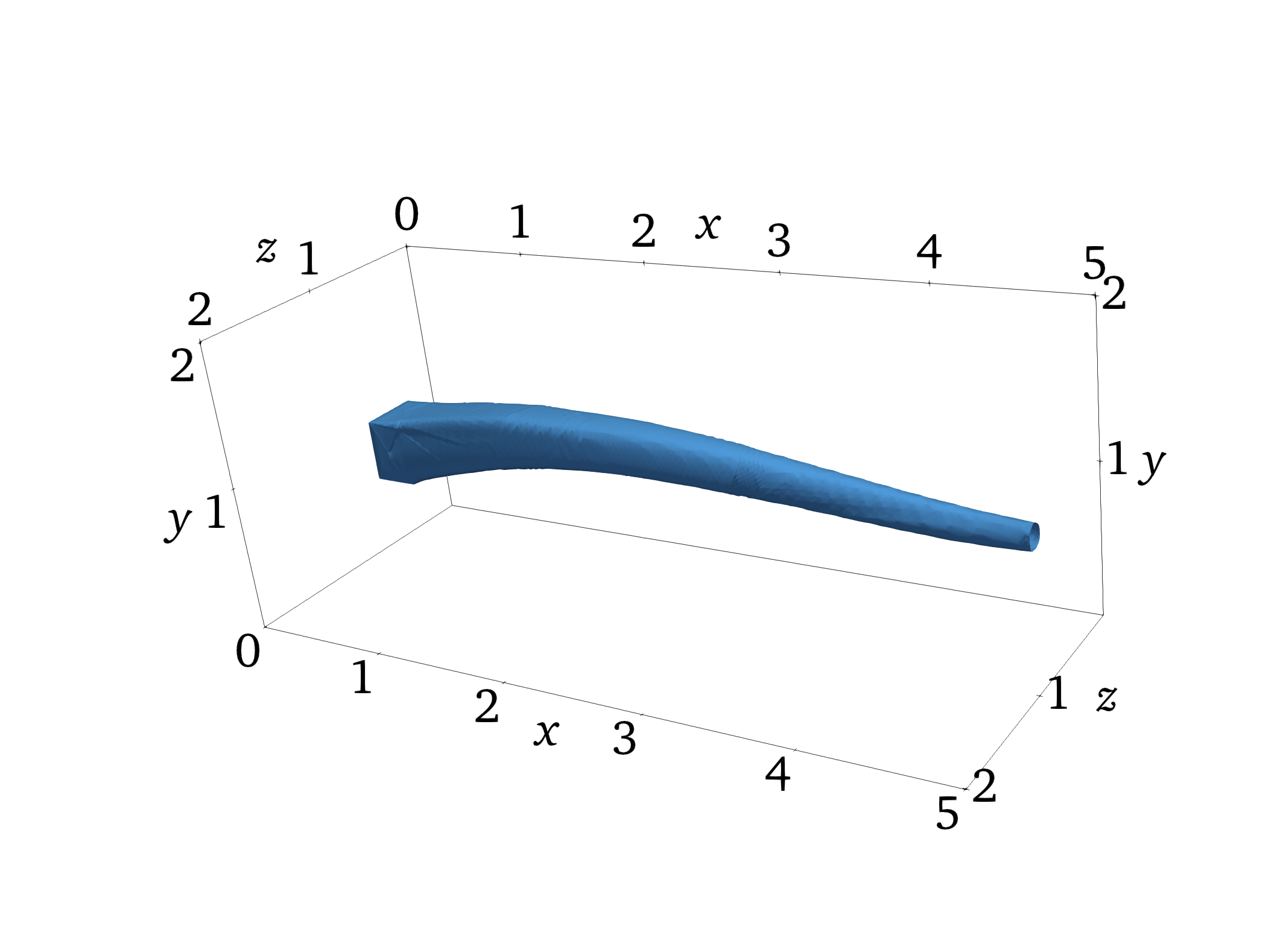}
\includegraphics[width=0.5\textwidth]{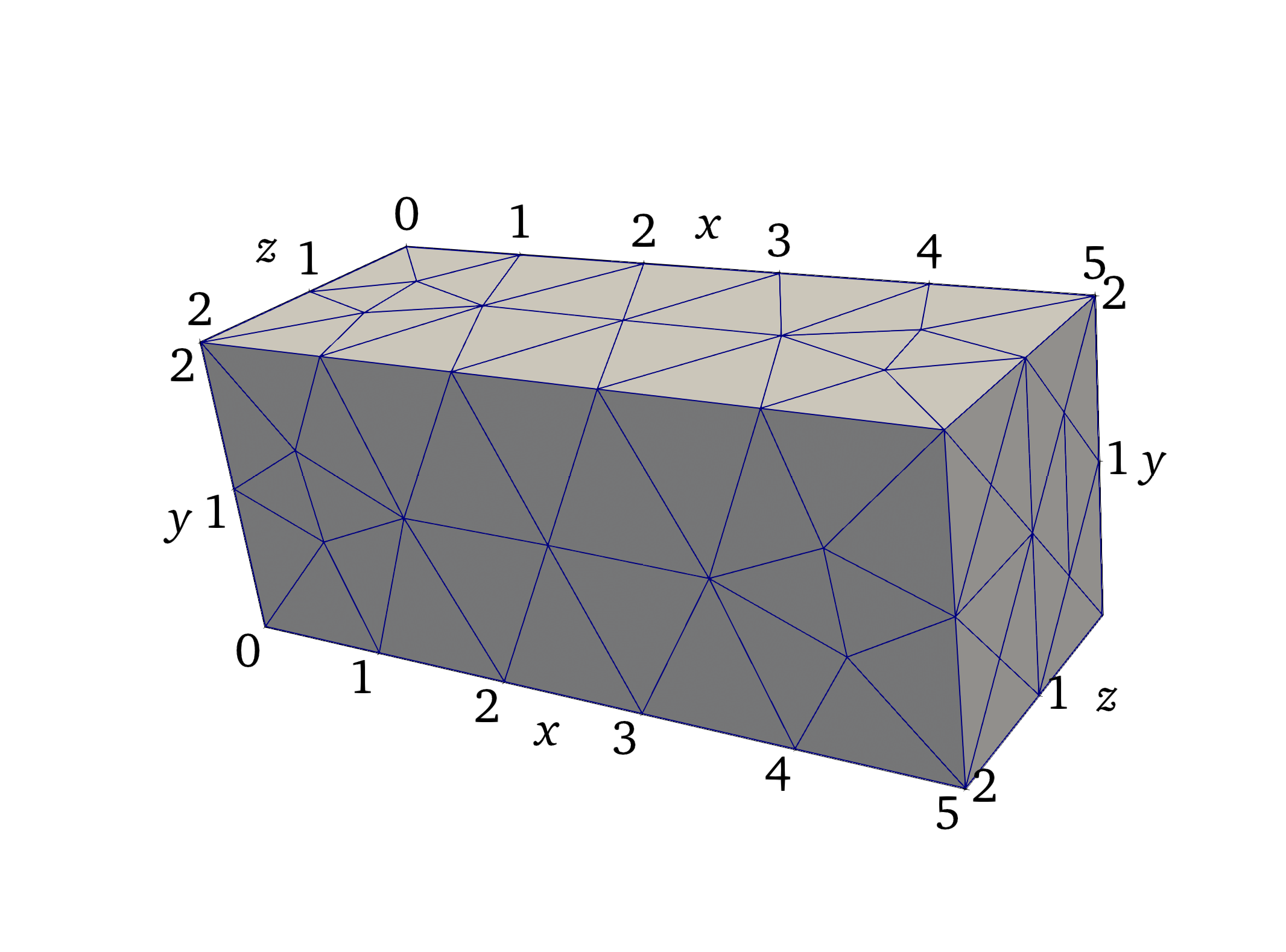}
\label{fig:bjkr18_sol}
\end{figure}

We consider ~\eqref{eq:cdr} with $\Omega = \Omega_1 \setminus \overline{\Omega}_2$, where $\Omega_1=(0,5)\times (0,2)\times (0,2)$ and $\Omega_2=(0.5,0.8)\times (0.8, 1.2) \times (0.8,1.2)$, $\bfb =(1, l(x), l(x))^T$ with $l(x)=(0.19x^3-1.42x^2-2.38x)/4$, $c = 0$, $f=0$ and the convection-dominated case of $\varepsilon = 10^{-3}$. An illustration of the solution is given in Figure~\ref{fig:bjkr18_sol}.

As for the boundary conditions, we set $\Gamma\nb := \lbrace 5 \rbrace \times (0,2) \times (0,2)$ and $\Gamma\db = \partial \Omega \setminus \Gamma\nb$, and we prescribe $g\nb = 0$ and $u\db = 0$ on $\partial \Omega_2$ and $u\db = 1$ on the remainder of $\Gamma\db$.

The region $\Omega$ is covered by unstructured tetrahedral grids corresponding to the refinement levels 3, 4 and 5. The coarsest one (in Figure~\ref{fig:bjkr18_sol}) was obtained by \textsc{Gmsh} \cite{GR09}. It consists of 226 tetrahedra. The cell diameters of our grids are in the ranges [0.0718, 0.2313], [0.0510, 0.1548] and [0.0180, 0.0578], respectively.
Finally, we consider $\mathtt{atol} = 10^{-14}$ and $\mathtt{tol} = 10^{-10}$.

Simulations were performed for all solvers and preconditioners listed in Tables~\ref{tab:solvers_iterative} and \ref{tab:pc} for all NP listed in item \ref{enum:parallelization} on page \pageref{enum:parallelization}. We note that certain combinations of solvers and preconditioners did not work; in particular, the \textsc{PETSc} solver with the preconditioner MG always either crashed or did not converge. Similarly, the \textsc{PETSc} solver always crashed when combining BCGS with Jac or SOR, or when combining LGMRES with SOR.

\subsubsection*{Discussion of Results}

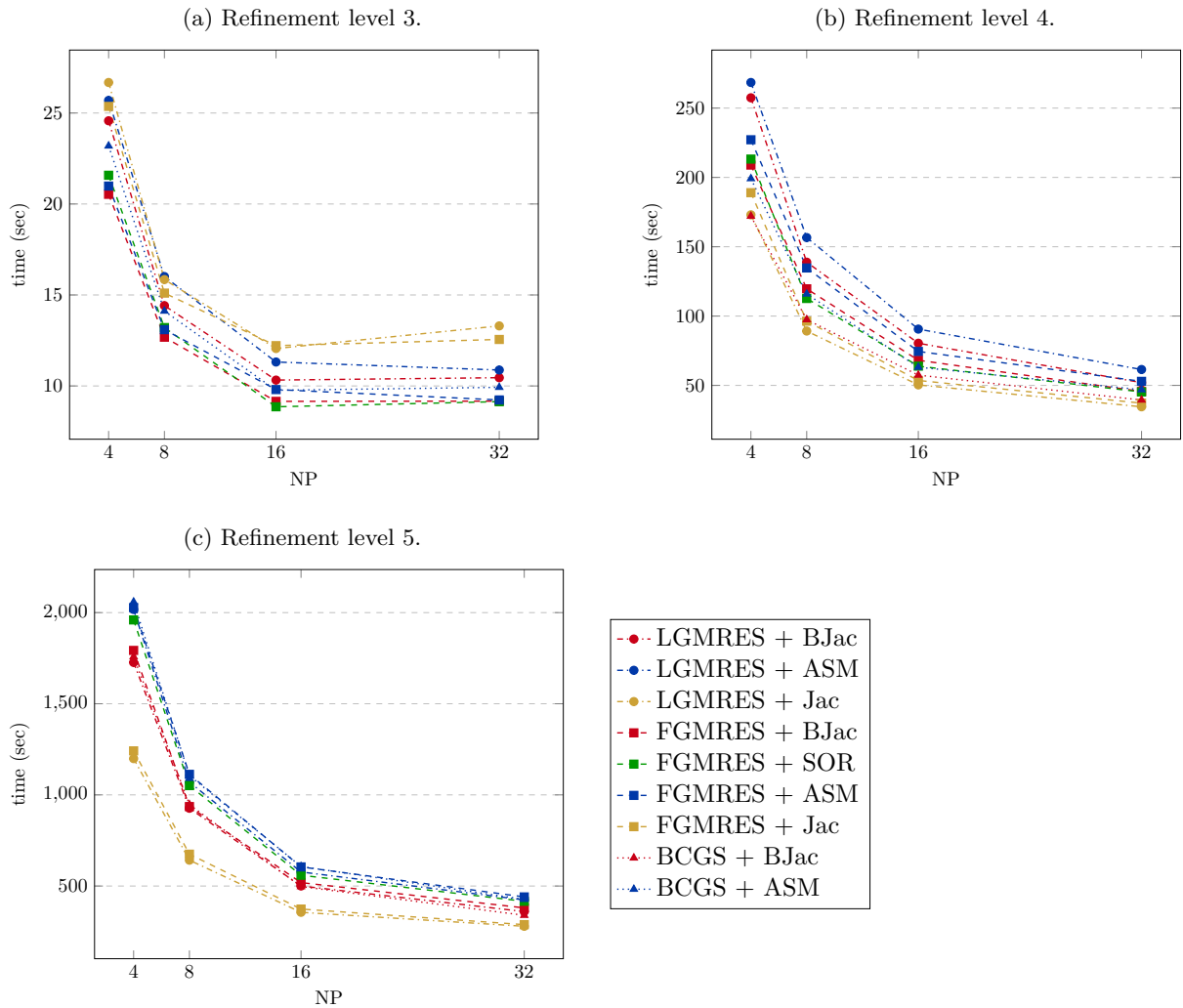
\begin{figure}[t!]
\caption{Example~\ref{ex:non_constant_convection2}: Efficiency of different solvers and preconditioners with respect to the number of processors for the LP limiter.}
\begin{subfigure}{0.49\textwidth}
\caption{Refinement level 3.}
\input{graphs/bjkr18/bjkr18_BJK17_3_figure3.tex}
\label{fig:bjk17_3}
\end{subfigure}\hspace*{1em}
\begin{subfigure}{0.49\textwidth}
\caption{Refinement level 4.}
\input{graphs/bjkr18/bjkr18_BJK17_4_figure3.tex}
\label{fig:bjk17_4}
\end{subfigure}
\par\bigskip
\begin{subfigure}{0.49\textwidth}
\caption{Refinement level 5.}
\input{graphs/bjkr18/bjkr18_BJK17_5_figure3.tex}
\label{fig:bjk17_5}
\end{subfigure}
\label{fig:bjk17}
\end{figure}

\begin{figure}[t!]
\caption{Example~\ref{ex:non_constant_convection2}: Efficiency of different solvers and preconditioners with respect to the number of processors for the MC limiter.}
\begin{subfigure}{0.49\textwidth}
\caption{Refinement level 3.}
\input{graphs/bjkr18/bjkr18_monolithic_3_figure3.tex}
\label{fig:monolithic_3}
\end{subfigure}\hspace*{1em}
\begin{subfigure}{0.49\textwidth}
\caption{Refinement level 4.}
\input{graphs/bjkr18/bjkr18_monolithic_4_figure3.tex}
\label{fig:monolithic_4}
\end{subfigure}
\par\bigskip
\begin{subfigure}{0.49\textwidth}
\caption{Refinement level 5.}
\input{graphs/bjkr18/bjkr18_monolithic_5_figure3.tex}
\label{fig:monolithic_5}
\end{subfigure}
\label{fig:monolithic}
\end{figure}
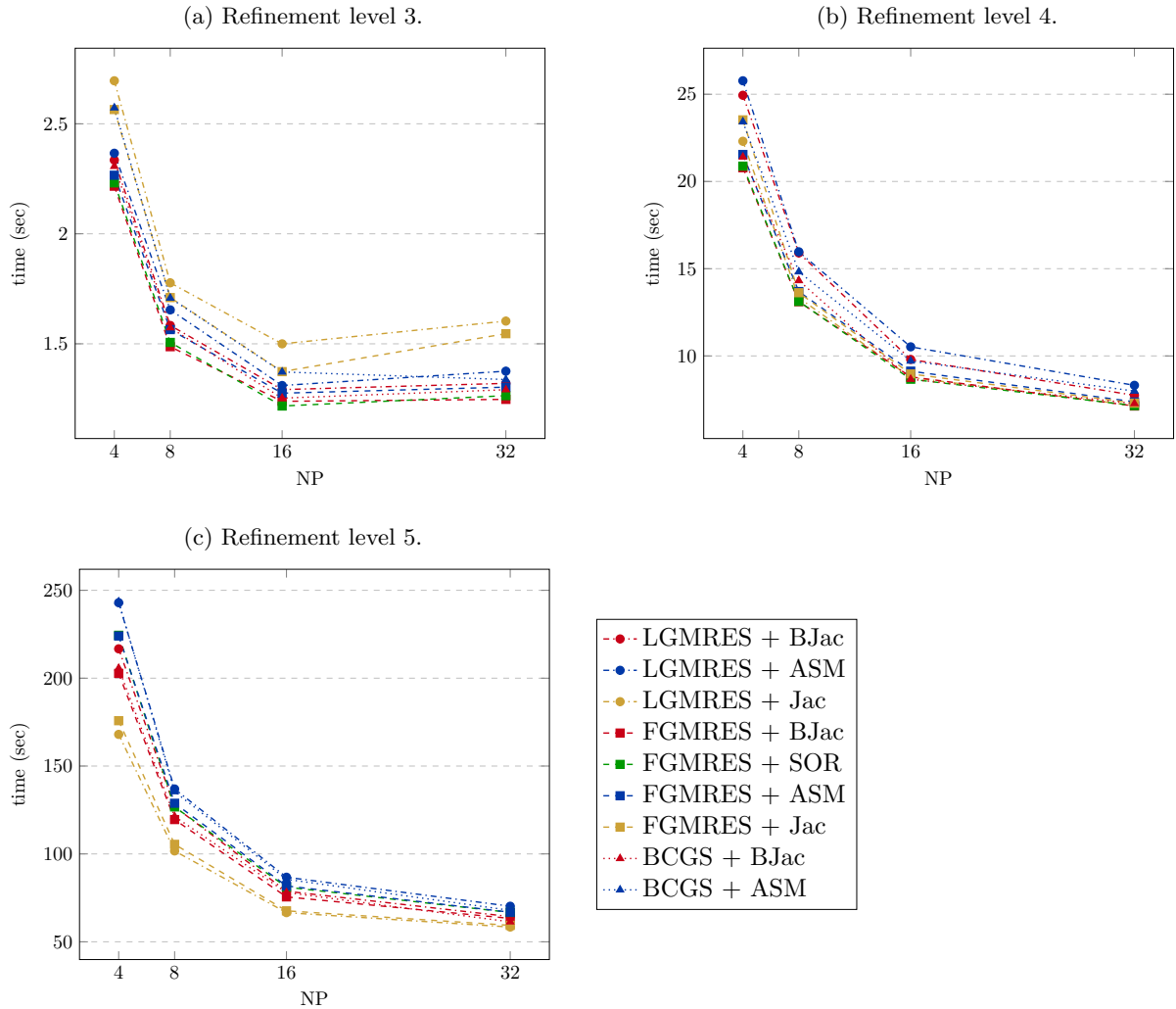

Results regarding the comparison of computing times for different limiters, refinement levels and number of processors (4, 8, 16 and 32) are compared in Figures~\ref{fig:bjk17}--\ref{fig:monolithic}, pages \pageref{fig:bjk17}--\pageref{fig:monolithic}.

\begin{figure}[t!]
\caption{Example~\ref{ex:non_constant_convection2}: Efficiency of different solvers and preconditioners with respect to the number of processors for the MU limiter.}
\begin{subfigure}{0.49\textwidth}
\caption{Refinement level 3.}
\input{graphs/bjkr18/bjkr18_kuzmin_3_figure3.tex}
\label{fig:kuzmin_3}
\end{subfigure}\hspace*{1em}
\begin{subfigure}{0.49\textwidth}
\caption{Refinement level 4.}
\input{graphs/bjkr18/bjkr18_kuzmin_4_figure3.tex}
\label{fig:kuzmin_4}
\end{subfigure}
\par\bigskip
\begin{subfigure}{0.49\textwidth}
\caption{Refinement level 5.}
\input{graphs/bjkr18/bjkr18_kuzmin_5_figure3.tex}
\label{fig:kuzmin_5}
\end{subfigure}
\label{fig:kuzmin}
\end{figure}
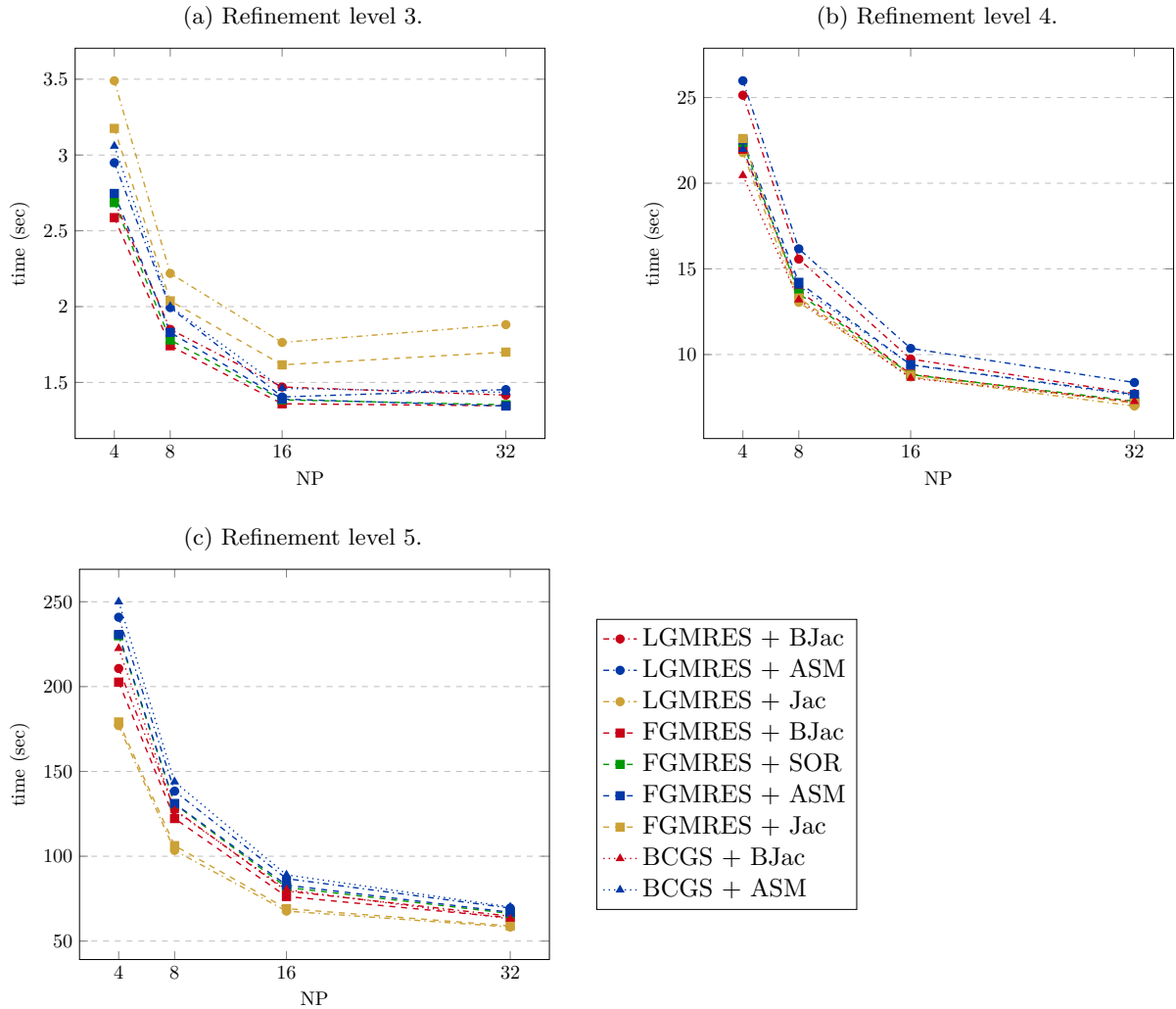

First, we believe that the shortening of the computing time by doubling the number of processors is acceptable. Note that for the lowest refinement level, the problem was obviously too small for 32 processors, and the computing time even increased sometimes.

Regarding the solvers, we note that FGMRES is clearly the most efficient. Further, BCGS seems to be slightly better than LGMRES.

As for the preconditioners, the best one is clearly Jac, but only when used for very fine grids: when coarsening the grid, it becomes less and less efficient.
The second best are BJac and SOR: Sometimes, BJac is better than SOR; sometimes SOR is the better one.

The worst preconditioners are clearly Jac (when used for coarse grids) and ASM.

The computing times of the scheme with the LP limiter are about an order of magnitude larger than those of the other runtimes. The same is true for the number of iterations (not shown here), which is the reason for these long runtimes. This seems to be the price to pay for better accuracy. This point is discussed in details in the next section.

In the discussion above, we compared only the parallel algorithms because we assumed that they would perform better than their sequential counterparts on the problem sizes considered. The correctness of this assumption is illustrated by Tables~\ref{tab:bjkr18_BJK17}--\ref{tab:bjkr18_monolithic}, where we have listed the sequential computing times and the parallel computing times for 4 processors for our lowest refinement level. The parallel implementations are clearly better.

\input{tables/bjkr18_BJK17_seqParTable3}
\input{tables/bjkr18_kuzmin_seqParTable3}
\input{tables/bjkr18_monolithic_seqParTable3}

\begin{figure}[t!]
\centering
\caption{Example~\ref{ex:non_constant_convection2}: $u(x,1,1)$ for $x \in [0,0.5] \cup [0.8,5]$ computed using different limiters.}
\input{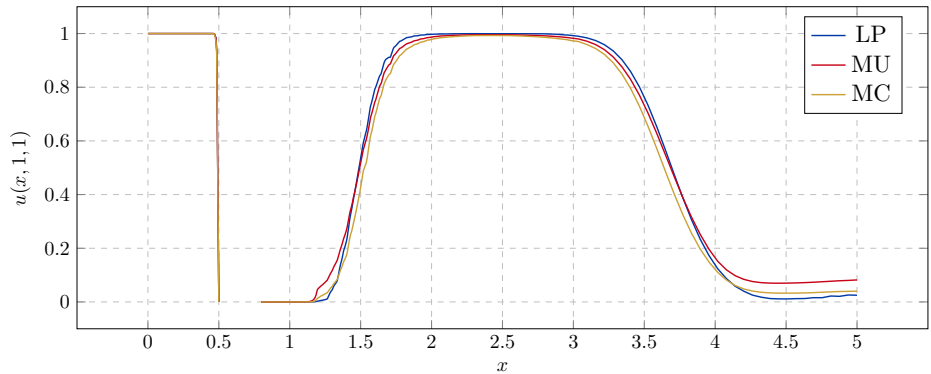}
\label{fig:cutLine}
\end{figure}

Finally, Figure~\ref{fig:cutLine} compares the solution values along the cut-line defined by $y = 1$ and $z = 1$. We can see that the solutions do not contain any overshoots or undershoots. They are also approximately the same for $x \leq 4.0$. However, they significantly differ for $x >4.0$.

%% file: graphs/bjkr18/bjkr18_BJK17_3_figure3.tex
\begin{tikzpicture}[scale=0.75]
\begin{axis}[
legend pos=north east, xlabel = NP, ylabel = {time (sec)},
legend style={nodes={scale=0.75, transform shape}},
xtick={4,8,16,32},
ymajorgrids=true,
grid style=dashed]
\addplot[color=crimson, line width = 0.25mm, mark=oplus*, mark options={solid}, dashdotted]
coordinates{
(4,24.5732)
(8,14.4143)
(16,10.3222)
(32,10.4538)
};
\addplot[color=royal_blue, line width = 0.25mm, mark=oplus*, mark options={solid}, dashdotted]
coordinates{
(4,25.692)
(8,16.011)
(16,11.3148)
(32,10.8752)
};
\addplot[color=gold, line width = 0.25mm, mark=oplus*, mark options={solid}, dashdotted]
coordinates{
(4,26.6773)
(8,15.8443)
(16,12.0569)
(32,13.2995)
};
\addplot[color=crimson, line width = 0.25mm, dashed, mark=square*, mark options={solid}]
coordinates{
(4,20.531)
(8,12.6791)
(16,9.15207)
(32,9.16565)
};
\addplot[color=dark_green, line width = 0.25mm, dashed, mark=square*, mark options={solid}]
coordinates{
(4,21.5794)
(8,13.207)
(16,8.85481)
(32,9.13789)
};
\addplot[color=royal_blue, line width = 0.25mm, dashed, mark=square*, mark options={solid}]
coordinates{
(4,20.9848)
(8,13.0925)
(16,9.80123)
(32,9.23042)
};
\addplot[color=gold, line width = 0.25mm, dashed, mark=square*, mark options={solid}]
coordinates{
(4,25.3628)
(8,15.1023)
(16,12.2141)
(32,12.5551)
};
\addplot[color=royal_blue, line width = 0.25mm, dotted, mark=triangle*, mark options={solid}]
coordinates{
(4,23.1895)
(8,14.1134)
(16,9.763)
(32,9.92374)
};
\end{axis}
\end{tikzpicture}

%% file: graphs/bjkr18/bjkr18_BJK17_4_figure3.tex
\begin{tikzpicture}[scale=0.75]
\begin{axis}[
legend pos=north east, xlabel = NP, ylabel = {time (sec)},
legend style={nodes={scale=0.75, transform shape}},
xtick={4,8,16,32},
ymajorgrids=true,
grid style=dashed]
\addplot[color=crimson, line width = 0.25mm, dashdotted, mark=oplus*, mark options={solid}]
coordinates{
(4,257.439)
(8,138.864)
(16,80.408)
(32,52.0096)
};
\addplot[color=royal_blue, line width = 0.25mm, dashdotted, mark=oplus*, mark options={solid}]
coordinates{
(4,268.464)
(8,156.609)
(16,90.5981)
(32,61.4537)
};
\addplot[color=gold, line width = 0.25mm, dashdotted, mark=oplus*, mark options={solid}]
coordinates{
(4,173.027)
(8,89.2047)
(16,50.3601)
(32,34.5735)
};
\addplot[color=crimson, line width = 0.25mm, dashed, mark=square*, mark options={solid}]
coordinates{
(4,208.895)
(8,119.655)
(16,68.1496)
(32,45.6066)
};
\addplot[color=dark_green, line width = 0.25mm, dashed, mark=square*, mark options={solid}]
coordinates{
(4,213.218)
(8,112.687)
(16,63.8823)
(32,45.2873)
};
\addplot[color=royal_blue, line width = 0.25mm, dashed, mark=square*, mark options={solid}]
coordinates{
(4,227.16)
(8,134.695)
(16,74.5031)
(32,52.9196)
};
\addplot[color=gold, line width = 0.25mm, dashed, mark=square*, mark options={solid}]
coordinates{
(4,188.989)
(8,96.3285)
(16,53.5438)
(32,37.1989)
};
\addplot[color=crimson, line width = 0.25mm, dotted, mark=triangle*, mark options={solid}]
coordinates{
(4,172.03)
(8,97.4096)
(16,57.4255)
(32,39.4523)
};
\addplot[color=royal_blue, line width = 0.25mm, dotted, mark=triangle*, mark options={solid}]
coordinates{
(4,199.116)
(8,115.908)
(16,63.1221)
(32,46.8518)
};
\end{axis}
\end{tikzpicture}

%% file: graphs/bjkr18/bjkr18_BJK17_5_figure3.tex
\begin{tikzpicture}[scale=0.75]
\begin{axis}[
legend pos=north east, xlabel = NP, ylabel = {time (sec)},
legend style={at={(1.1, 0.5)},anchor=west, nodes={scale=1.3, transform shape}},
legend cell align={left},
xtick={4,8,16,32},
ymajorgrids=true,
grid style=dashed]
\addplot[color=crimson, line width = 0.25mm, dashdotted, mark=oplus*, mark options={solid}]
coordinates{
(4,1726.62)
(8,926.596)
(16,500.266)
(32,359.856)
};
\addlegendentry{LGMRES + BJac}
\addplot[color=royal_blue, line width = 0.25mm, dashdotted, mark=oplus*, mark options={solid}]
coordinates{
(4,2017.3)
(8,1070.17)
(16,577.206)
(32,423.765)
};
\addlegendentry{LGMRES + ASM}
\addplot[color=gold, line width = 0.25mm, dashdotted, mark=oplus*, mark options={solid}]
coordinates{
(4,1199.06)
(8,641.748)
(16,355.777)
(32,279.201)
};
\addlegendentry{LGMRES + Jac}
\addplot[color=crimson, line width = 0.25mm, dashed, mark=square*, mark options={solid}]
coordinates{
(4,1792.48)
(8,935.803)
(16,518.754)
(32,380.023)
};
\addlegendentry{FGMRES + BJac}
\addplot[color=dark_green, line width = 0.25mm, dashed, mark=square*, mark options={solid}]
coordinates{
(4,1959.79)
(8,1050.98)
(16,559.213)
(32,416.236)
};
\addlegendentry{FGMRES + SOR}
\addplot[color=royal_blue, line width = 0.25mm, dashed, mark=square*, mark options={solid}]
coordinates{
(4,2027.35)
(8,1112.68)
(16,604.917)
(32,440.272)
};
\addlegendentry{FGMRES + ASM}
\addplot[color=gold, line width = 0.25mm, dashed, mark=square*, mark options={solid}]
coordinates{
(4,1241.44)
(8,674.738)
(16,373.71)
(32,288.23)
};
\addlegendentry{FGMRES + Jac}
\addplot[color=crimson, line width = 0.25mm, dotted, mark=triangle*, mark options={solid}]
coordinates{
(4,1758.23)
(8,946.745)
(16,496.997)
(32,339.443)
};
\addlegendentry{BCGS + BJac}
\addplot[color=royal_blue, line width = 0.25mm, dotted, mark=triangle*, mark options={solid}]
coordinates{
(4,2058.3)
(8,1105.33)
(16,607.13)
(32,426.811)
};
\addlegendentry{BCGS + ASM}
\end{axis}
\end{tikzpicture}

%% file: graphs/bjkr18/bjkr18_kuzmin_3_figure3.tex
\begin{tikzpicture}[scale=0.75]
\begin{axis}[
legend pos=north east, xlabel = NP, ylabel = {time (sec)},
legend style={nodes={scale=0.75, transform shape}},
xtick={4,8,16,32},
ymajorgrids=true,
grid style=dashed]
\addplot[color=crimson, line width = 0.25mm, dashdotted, mark=oplus*, mark options={solid}]
coordinates{
(4,2.69996)
(8,1.84984)
(16,1.4698)
(32,1.41634)
};
\addplot[color=royal_blue, line width = 0.25mm, dashdotted, mark=oplus*, mark options={solid}]
coordinates{
(4,2.94906)
(8,1.9936)
(16,1.40407)
(32,1.45261)
};
\addplot[color=gold, line width = 0.25mm, dashdotted, mark=oplus*, mark options={solid}]
coordinates{
(4,3.48971)
(8,2.2189)
(16,1.7637)
(32,1.88102)
};
\addplot[color=crimson, line width = 0.25mm, dashed, mark=square*, mark options={solid}]
coordinates{
(4,2.58716)
(8,1.74166)
(16,1.35901)
(32,1.34595)
};
\addplot[color=dark_green, line width = 0.25mm, dashed, mark=square*, mark options={solid}]
coordinates{
(4,2.68637)
(8,1.77978)
(16,1.38504)
(32,1.35261)
};
\addplot[color=royal_blue, line width = 0.25mm, dashed, mark=square*, mark options={solid}]
coordinates{
(4,2.74627)
(8,1.82986)
(16,1.39002)
(32,1.34433)
};
\addplot[color=gold, line width = 0.25mm, dashed, mark=square*, mark options={solid}]
coordinates{
(4,3.17527)
(8,2.03803)
(16,1.61587)
(32,1.69991)
};
\addplot[color=royal_blue, line width = 0.25mm, dotted, mark=triangle*, mark options={solid}]
coordinates{
(4,3.05882)
(8,2.00137)
(16,1.46077)
(32,1.43426)
};
\end{axis}
\end{tikzpicture}

%% file: graphs/bjkr18/bjkr18_kuzmin_4_figure3.tex
\begin{tikzpicture}[scale=0.75]
\begin{axis}[
legend pos=north east, xlabel = NP, ylabel = {time (sec)},
legend style={nodes={scale=0.75, transform shape}},
xtick={4,8,16,32},
ymajorgrids=true,
grid style=dashed]
\addplot[color=crimson, line width = 0.25mm, dashdotted, mark=oplus*, mark options={solid}]
coordinates{
(4,25.1383)
(8,15.5702)
(16,9.73921)
(32,7.7165)
};
\addplot[color=royal_blue, line width = 0.25mm, dashdotted, mark=oplus*, mark options={solid}]
coordinates{
(4,25.9799)
(8,16.1714)
(16,10.3589)
(32,8.36771)
};
\addplot[color=gold, line width = 0.25mm, dashdotted, mark=oplus*, mark options={solid}]
coordinates{
(4,21.789)
(8,13.0473)
(16,8.69222)
(32,6.99593)
};
\addplot[color=crimson, line width = 0.25mm, dashed, mark=square*, mark options={solid}]
coordinates{
(4,21.9914)
(8,13.7305)
(16,8.86198)
(32,7.18933)
};
\addplot[color=dark_green, line width = 0.25mm, dashed, mark=square*, mark options={solid}]
coordinates{
(4,22.3785)
(8,13.5658)
(16,8.82796)
(32,7.27601)
};
\addplot[color=royal_blue, line width = 0.25mm, dashed, mark=square*, mark options={solid}]
coordinates{
(4,22.524)
(8,14.2161)
(16,9.40656)
(32,7.67095)
};
\addplot[color=gold, line width = 0.25mm, dashed, mark=square*, mark options={solid}]
coordinates{
(4,22.6062)
(8,13.301)
(16,8.79463)
(32,7.16608)
};
\addplot[color=crimson, line width = 0.25mm, dotted, mark=triangle*, mark options={solid}]
coordinates{
(4,20.4555)
(8,13.2001)
(16,8.62279)
(32,7.26577)
};
\addplot[color=royal_blue, line width = 0.25mm, dotted, mark=triangle*, mark options={solid}]
coordinates{
(4,21.9848)
(8,14.0069)
(16,9.40059)
(32,7.61933)
};
\end{axis}
\end{tikzpicture}

%% file: graphs/bjkr18/bjkr18_kuzmin_5_figure3.tex
\begin{tikzpicture}[scale=0.75]
\begin{axis}[
legend pos=north east, xlabel = NP, ylabel = {time (sec)},
legend style={at={(1.1, 0.5)},anchor=west, nodes={scale=1.3, transform shape}},
legend cell align={left},
xtick={4,8,16,32},
ymajorgrids=true,
grid style=dashed]
\addplot[color=crimson, line width = 0.25mm, dashdotted, mark=oplus*, mark options={solid}]
coordinates{
(4,210.65)
(8,127.539)
(16,79.5323)
(32,64.6589)
};
\addlegendentry{LGMRES + BJac}
\addplot[color=royal_blue, line width = 0.25mm, dashdotted, mark=oplus*, mark options={solid}]
coordinates{
(4,240.864)
(8,138.368)
(16,86.7985)
(32,69.4481)
};
\addlegendentry{LGMRES + ASM}
\addplot[color=gold, line width = 0.25mm, dashdotted, mark=oplus*, mark options={solid}]
coordinates{
(4,176.984)
(8,103.427)
(16,67.665)
(32,58.2836)
};
\addlegendentry{LGMRES + Jac}
\addplot[color=crimson, line width = 0.25mm, dashed, mark=square*, mark options={solid}]
coordinates{
(4,202.62)
(8,122.232)
(16,76.27)
(32,63.7546)
};
\addlegendentry{FGMRES + BJac}
\addplot[color=dark_green, line width = 0.25mm, dashed, mark=square*, mark options={solid}]
coordinates{
(4,229.953)
(8,130.805)
(16,81.7068)
(32,66.2457)
};
\addlegendentry{FGMRES + SOR}
\addplot[color=royal_blue, line width = 0.25mm, dashed, mark=square*, mark options={solid}]
coordinates{
(4,230.714)
(8,131.096)
(16,83.0036)
(32,66.943)
};
\addlegendentry{FGMRES + ASM}
\addplot[color=gold, line width = 0.25mm, dashed, mark=square*, mark options={solid}]
coordinates{
(4,179.168)
(8,106.264)
(16,69.1173)
(32,58.9391)
};
\addlegendentry{FGMRES + Jac}
\addplot[color=crimson, line width = 0.25mm, dotted, mark=triangle*, mark options={solid}]
coordinates{
(4,222.509)
(8,127.006)
(16,80.3232)
(32,63.0596)
};
\addlegendentry{BCGS + BJac}
\addplot[color=royal_blue, line width = 0.25mm, dotted, mark=triangle*, mark options={solid}]
coordinates{
(4,249.981)
(8,143.91)
(16,89.0346)
(32,70.096)
};
\addlegendentry{BCGS + ASM}
\end{axis}
\end{tikzpicture}

%% file: graphs/bjkr18/bjkr18_monolithic_3_figure3.tex
\begin{tikzpicture}[scale=0.75]
\begin{axis}[
legend pos=north east, xlabel = NP, ylabel = {time (sec)},
legend style={nodes={scale=0.75, transform shape}},
xtick={4,8,16,32},
ymajorgrids=true,
grid style=dashed]
\addplot[color=crimson, line width = 0.25mm, dashdotted, mark=oplus*, mark options={solid}]
coordinates{
(4,2.33499)
(8,1.58348)
(16,1.29178)
(32,1.32022)
};
\addplot[color=royal_blue, line width = 0.25mm, dashdotted, mark=oplus*, mark options={solid}]
coordinates{
(4,2.36623)
(8,1.65414)
(16,1.31049)
(32,1.37558)
};
\addplot[color=gold, line width = 0.25mm, dashdotted, mark=oplus*, mark options={solid}]
coordinates{
(4,2.69557)
(8,1.7771)
(16,1.4998)
(32,1.60321)
};
\addplot[color=crimson, line width = 0.25mm, dashed, mark=square*, mark options={solid}]
coordinates{
(4,2.2164)
(8,1.48675)
(16,1.23853)
(32,1.24738)
};
\addplot[color=dark_green, line width = 0.25mm, dashed, mark=square*, mark options={solid}]
coordinates{
(4,2.23136)
(8,1.50634)
(16,1.2172)
(32,1.26424)
};
\addplot[color=royal_blue, line width = 0.25mm, dashed, mark=square*, mark options={solid}]
coordinates{
(4,2.26615)
(8,1.56459)
(16,1.27554)
(32,1.30291)
};
\addplot[color=gold, line width = 0.25mm, dashed, mark=square*, mark options={solid}]
coordinates{
(4,2.56485)
(8,1.71011)
(16,1.37347)
(32,1.54596)
};
\addplot[color=crimson, line width = 0.25mm, dotted, mark=triangle*, mark options={solid}]
coordinates{
(4,2.30653)
(8,1.57295)
(16,1.25307)
(32,1.29128)
};
\addplot[color=royal_blue, line width = 0.25mm, dotted, mark=triangle*, mark options={solid}]
coordinates{
(4,2.57253)
(8,1.7079)
(16,1.372)
(32,1.33761)
};
\end{axis}
\end{tikzpicture}

%% file: graphs/bjkr18/bjkr18_monolithic_4_figure3.tex
\begin{tikzpicture}[scale=0.75]
\begin{axis}[
legend pos=north east, xlabel = NP, ylabel = {time (sec)},
legend style={nodes={scale=0.75, transform shape}},
xtick={4,8,16,32},
ymajorgrids=true,
grid style=dashed]
\addplot[color=crimson, line width = 0.25mm, dashdotted, mark=oplus*, mark options={solid}]
coordinates{
(4,24.9415)
(8,15.9037)
(16,9.79745)
(32,7.74209)
};
\addplot[color=royal_blue, line width = 0.25mm, dashdotted, mark=oplus*, mark options={solid}]
coordinates{
(4,25.7717)
(8,15.9683)
(16,10.5203)
(32,8.31448)
};
\addplot[color=gold, line width = 0.25mm, dashdotted, mark=oplus*, mark options={solid}]
coordinates{
(4,22.3042)
(8,13.3296)
(16,8.76624)
(32,7.15362)
};
\addplot[color=crimson, line width = 0.25mm, dashed, mark=square*, mark options={solid}]
coordinates{
(4,20.7881)
(8,13.1075)
(16,8.79992)
(32,7.12805)
};
\addplot[color=dark_green, line width = 0.25mm, dashed, mark=square*, mark options={solid}]
coordinates{
(4,20.8773)
(8,13.1262)
(16,8.66416)
(32,7.14922)
};
\addplot[color=royal_blue, line width = 0.25mm, dashed, mark=square*, mark options={solid}]
coordinates{
(4,21.5439)
(8,13.689)
(16,9.12674)
(32,7.36465)
};
\addplot[color=gold, line width = 0.25mm, dashed, mark=square*, mark options={solid}]
coordinates{
(4,23.5303)
(8,13.623)
(16,8.96811)
(32,7.2928)
};
\addplot[color=crimson, line width = 0.25mm, dotted, mark=triangle*, mark options={solid}]
coordinates{
(4,21.3997)
(8,14.3182)
(16,8.66379)
(32,7.26977)
};
\addplot[color=royal_blue, line width = 0.25mm, dotted, mark=triangle*, mark options={solid}]
coordinates{
(4,23.4192)
(8,14.8237)
(16,9.70626)
(32,7.98376)
};
\end{axis}
\end{tikzpicture}

%% file: graphs/bjkr18/bjkr18_monolithic_5_figure3.tex
\begin{tikzpicture}[scale=0.75]
\begin{axis}[
legend pos=north east, xlabel = NP, ylabel = {time (sec)},
legend style={at={(1.1, 0.5)},anchor=west, nodes={scale=1.3, transform shape}},
legend cell align={left},
xtick={4,8,16,32},
ymajorgrids=true,
grid style=dashed]
\addplot[color=crimson, line width = 0.25mm, dashdotted, mark=oplus*, mark options={solid}]
coordinates{
(4,216.705)
(8,126.858)
(16,78.7698)
(32,64.3283)
};
\addlegendentry{LGMRES + BJac}
\addplot[color=royal_blue, line width = 0.25mm, dashdotted, mark=oplus*, mark options={solid}]
coordinates{
(4,242.954)
(8,136.959)
(16,86.7429)
(32,70.2382)
};
\addlegendentry{LGMRES + ASM}
\addplot[color=gold, line width = 0.25mm, dashdotted, mark=oplus*, mark options={solid}]
coordinates{
(4,168.016)
(8,101.697)
(16,66.5842)
(32,58.3972)
};
\addlegendentry{LGMRES + Jac}
\addplot[color=crimson, line width = 0.25mm, dashed, mark=square*, mark options={solid}]
coordinates{
(4,202.732)
(8,119.548)
(16,75.5414)
(32,63.2999)
};
\addlegendentry{FGMRES + BJac}
\addplot[color=dark_green, line width = 0.25mm, dashed, mark=square*, mark options={solid}]
coordinates{
(4,224.42)
(8,126.75)
(16,81.1014)
(32,66.6922)
};
\addlegendentry{FGMRES + SOR}
\addplot[color=royal_blue, line width = 0.25mm, dashed, mark=square*, mark options={solid}]
coordinates{
(4,224.009)
(8,128.832)
(16,81.786)
(32,66.831)
};
\addlegendentry{FGMRES + ASM}
\addplot[color=gold, line width = 0.25mm, dashed, mark=square*, mark options={solid}]
coordinates{
(4,175.803)
(8,105.449)
(16,67.6123)
(32,59.3598)
};
\addlegendentry{FGMRES + Jac}
\addplot[color=crimson, line width = 0.25mm, dotted, mark=triangle*, mark options={solid}]
coordinates{
(4,205.636)
(8,121.309)
(16,78.0358)
(32,61.3608)
};
\addlegendentry{BCGS + BJac}
\addplot[color=royal_blue, line width = 0.25mm, dotted, mark=triangle*, mark options={solid}]
coordinates{
(4,243.477)
(8,135.422)
(16,85.5431)
(32,67.9064)
};
\addlegendentry{BCGS + ASM}
\end{axis}
\end{tikzpicture}

%% file: tables/bjkr18_BJK17_seqParTable3.tex
\begin{table}[t!]
\centering
\caption{Example~\ref{ex:non_constant_convection2}: Computing times in seconds of the solver with the LP limiter for refinement level 3.}
\begin{tabular}{c c c c c c c c}
\multicolumn{2}{c}{} & \multicolumn{2}{c}{\textbf{time (sec)}}&\multicolumn{2}{c}{} &\multicolumn{2}{c}{\textbf{time (sec)}}\\
\textbf{Solver} & \textbf{PC} & NP = 1 & NP = 4 & \textbf{Solver} & \textbf{PC} & NP = 1 & NP = 4\\ \hline
\multirow{4}{*}{LGMRES} & BJac & 74.6 & 24.6 & \multirow{3}{*}{BCGS} & BJac & 63.1 & 21.7 \\ \cline{2-4} \cline{6-8}
 & SOR  & 79.6 & ---  & & SOR    & 66.8 & --- \\ \cline{2-4} \cline{6-8}
 & ASM  & 79.8 & 25.7 & & ASM    & 68.4 & 23.2 \\ \cline{2-4} \cline{5-8}
 & Jac  & 91.0 & 26.7 & &     &     &     \\ \hline
\multirow{4}{*}{FGMRES} & BJac & 60.8 & 20.5 & & & & \\ \cline{2-4} \cline{5-8}
 & SOR  & 62.9 & 21.6 & & & &\\ \cline{2-4} \cline{5-8}
 & ASM  & 65.3 & 21.0 & & & &\\ \cline{2-4} \cline{5-8}
 & Jac  & 84.9 & 25.4 & & & &\\ 
\end{tabular}
\label{tab:bjkr18_BJK17}
\end{table}

%% file: tables/bjkr18_kuzmin_seqParTable3.tex
\begin{table}[t!]
\centering
\caption{Example~\ref{ex:non_constant_convection2}: Computing times in seconds of the solver with the MU limiter for refinement level 3.}
\begin{tabular}{c c c c c c c c}
\multicolumn{2}{c}{} & \multicolumn{2}{c}{\textbf{time (sec)}}&\multicolumn{2}{c}{} &\multicolumn{2}{c}{\textbf{time (sec)}}\\
\textbf{Solver} & \textbf{PC} & NP = 1 & NP = 4 & \textbf{Solver} & \textbf{PC} & NP = 1 & NP = 4\\ \hline
\multirow{4}{*}{LGMRES} & BJac & 8.2 & 2.7 & \multirow{3}{*}{BCGS} & BJac & 7.9 & 2.9 \\ \cline{2-4} \cline{6-8}
 & SOR & 8.5 & ---  & & SOR & 8.1 & --- \\ \cline{2-4} \cline{6-8}
 & ASM & 8.7 & 2.9  & & ASM & 8.3 & 3.1 \\ \cline{2-4} \cline{5-8}
 & Jac & 11.1 & 3.5 & &     &     &     \\ \hline
\multirow{4}{*}{FGMRES} & BJac & 7.6 & 2.6 & & & & \\ \cline{2-4} \cline{5-8}
 & SOR & 7.8 & 2.7  & & & &\\ \cline{2-4} \cline{5-8}
 & ASM & 8.0 & 2.7  & & & &\\ \cline{2-4} \cline{5-8}
 & Jac & 10.0 & 3.2 & & & &\\
\end{tabular}
\label{tab:bjkr18_kuzmin}
\end{table}

%% file: tables/bjkr18_monolithic_seqParTable3.tex
\begin{table}[t!]
\centering
\caption{Example~\ref{ex:non_constant_convection2}: Computing times in seconds of the solver with the MC limiter for refinement level 3.}
\begin{tabular}{c c c c c c c c}
\multicolumn{2}{c}{} & \multicolumn{2}{c}{\textbf{time (sec)}}&\multicolumn{2}{c}{} &\multicolumn{2}{c}{\textbf{time (sec)}}\\
\textbf{Solver} & \textbf{PC} & NP = 1 & NP = 4 & \textbf{Solver} & \textbf{PC} & NP = 1 & NP = 4\\ \hline
\multirow{4}{*}{LGMRES} & BJac & 6.2 & 2.3 & \multirow{3}{*}{BCGS} & BJac & 5.8 & 2.3 \\ \cline{2-4} \cline{6-8}
 & SOR  & 6.4 & --- & & SOR    & 6.0 & --- \\ \cline{2-4} \cline{6-8}
 & ASM  & 6.6 & 2.4 & & ASM    & 6.2 & 2.6 \\ \cline{2-4} \cline{5-8}
 & Jac  & 8.1 & 2.7 & &     &     &     \\ \hline
\multirow{4}{*}{FGMRES} & BJac & 5.7 & 2.2 & & & & \\ \cline{2-4} \cline{5-8} 
 & SOR & 5.8 & 2.2 & & & &\\ \cline{2-4} \cline{5-8}
 & ASM & 6.0 & 2.3 & & & &\\ \cline{2-4} \cline{5-8}
 & Jac & 7.5 & 2.6 & & & &\\
\end{tabular}
\label{tab:bjkr18_monolithic}
\end{table}

%% file: circular_convection.tex
This example is an extension of the 2D example with the same name from \cite{Ku20}. We use it to compare the accuracy of the solutions computed with different limiters.

\subsubsection*{Description of Problem}

We consider equation \eqref{eq:cdr} with $\Omega = (0,1)^3$, $\varepsilon = 0$, $\bfb = (y, -x, 0)^\text{T}$, $c = 0$, and $f=0$.

\begin{figure}[t!]
\centering
\caption{Example~\ref{ex:circular_convection}: Solution for the refinement level 7 computed by the scheme with the MC limiter.}
\includegraphics[width=0.5\textwidth, trim = 3cm 0cm 8cm 9cm, clip = true]
{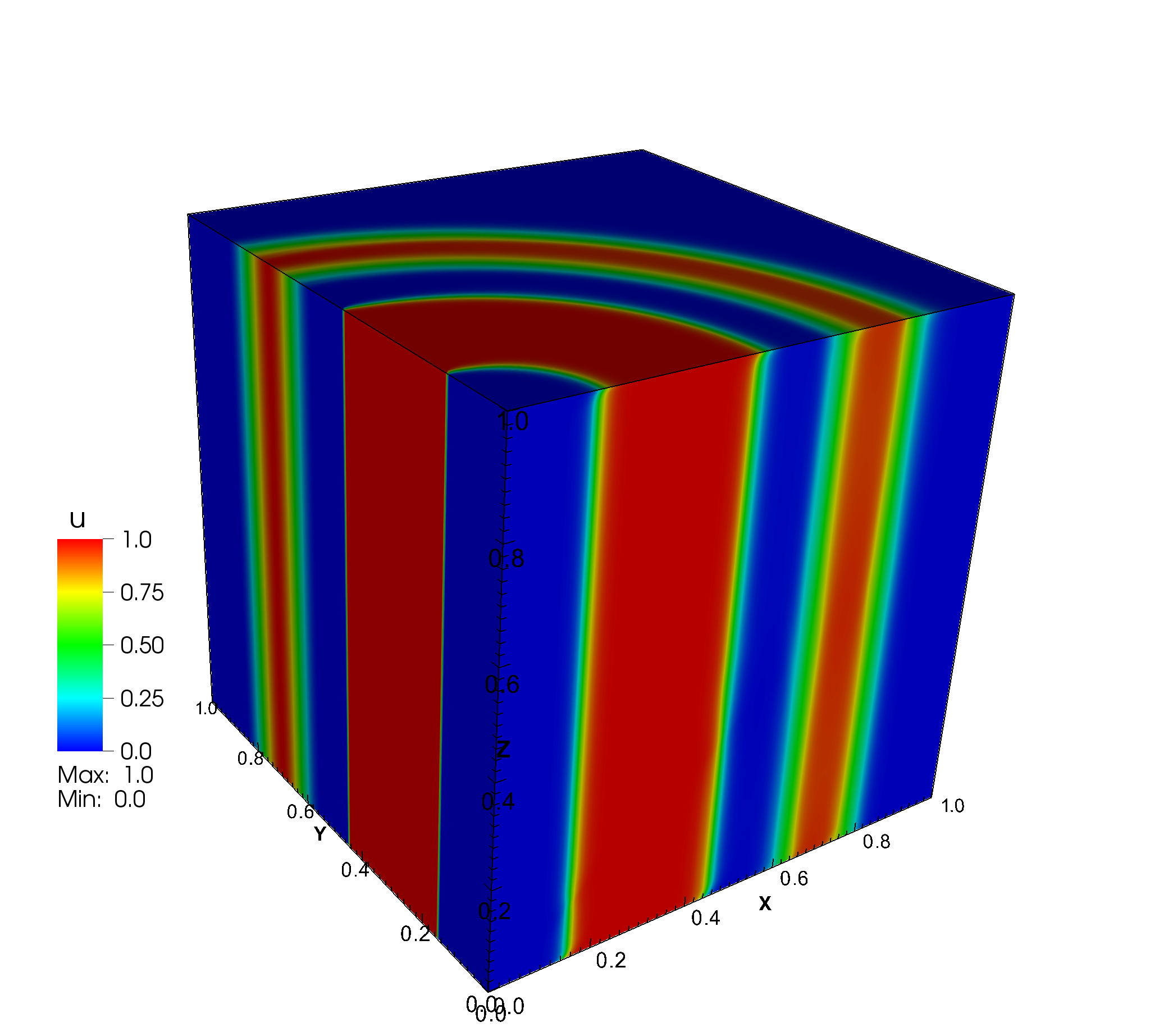}
\label{fig:circular_convection_mcl}
\end{figure}

As for the boundary conditions, $\Gamma_D$ is the union of the faces with $y=1$, $x=0$ and $x=1$, and $\Gamma_N = \Gamma \setminus \Gamma_D$, where $g\nb=0$. The values of $u\db$ correspond to the exact solution
\begin{equation}
u(x,y,z) = \begin{cases} 1 \ \ &\text{if} \ \ 0.15 \leq r(x,y) \leq 0.45, \\
\cos^2 \left( 10 \pi \frac{r(x,y) - 0.7}{3} \right) \ \ &\text{if} \ \ 0.55 \leq r(x,y) \leq 0.85, \\
0 \ \ &\text{otherwise},
 \end{cases}
\end{equation}
where $r(x,y) = \sqrt{x^2 + y^2}$. An approximation of this solution is in Figure \ref{fig:circular_convection_mcl}.

The domain $\Omega$ is covered by tetrahedral grids corresponding to the refinement levels 5--7, with the initial grid consisting of 6 tetrahedra. The maximum cell diameters of these grids are approximately $0.0765, 0.0383$, and $0.0191$, respectively.

Since we wanted to obtain the best approximation of the solution for a given combination of limiter and refinement level, our stopping criterion for the nonlinear solver was
\begin{equation}
\label{eq:circular_convection_stopping_criterion}
|R_\text{new} - R_\text{old}| / R_\text{new} < \mathtt{tol}_2 \ \ \text{and} \ \ u(\bfx) \geq -10^{-16} \ \forall \bfx \in \overline{\Omega} ,
\end{equation}
where $R_{\mathrm{new}}$ and $R_{\mathrm{old}}$ stand for the Euclidean norms of the new and old residues, respectively. That is, the solver should stop when the solution stops improving. We set $\mathtt{tol}_2 = 10^{-6}$ and $\mathtt{atol} = 10^{-14}$.

\subsubsection*{Discussion of Results}

Table~\ref{tab:circular_convection_errors}, page~\pageref{tab:circular_convection_errors} shows the resulting errors of the numerical solution, measured in $\left\Arrowvert \cdot \right\Arrowvert_1$ and $\left\Arrowvert \cdot \right\Arrowvert_2$, which denote the standard norms in  $L^1(\Omega)$ and $L^2(\Omega)$, respectively. Note that we do not specify which combination of linear systems solver and preconditioner we use.
The reason for this is that we tested all the solvers and preconditioners listed in Tables~\ref{tab:solvers_iterative} and \ref{tab:pc}, and, as expected, they all gave approximately the same results (with the exception of some combinations, for which the \textsc{PETSc} solver always crashed). In addition, we only ran our solvers in parallel.

The errors clearly indicate convergence. However, the $L^1$-error decreases much faster than the $L^2$-error. For each refinement level, the results obtained with the LP limiter are clearly the best. The results obtained with the MC and MU limiters are similar.

Due to \eqref{eq:circular_convection_stopping_criterion}, none of the schemes produced undershoots or overshoots that were larger than the machine precision.

\input{tables/circular_convection_L1_L2_precTable.tex}

%% file: tables/circular_convection_L1_L2_precTable.tex
\begin{table}[t!]
\centering
\caption{Example~\ref{ex:circular_convection}: Errors measured in $\left\Arrowvert \cdot \right\Arrowvert_1$ and $\left\Arrowvert \cdot \right\Arrowvert_2$.}
\begin{subtable}{\textwidth}
\centering
\caption{Error in $\left\Arrowvert \cdot \right\Arrowvert_1$.}
\begin{tabular}{c c c c}
\textbf{Level} & \textbf{MC} & \textbf{LP} & \textbf{MU} \\ \hline
5 & 8.11e-02 & 4.14e-02 & 9.24e-02 \\ 
6 & 3.39e-02 & 1.86e-02 & 3.62e-02 \\
7 & 1.62e-02 & 8.31e-03 & 1.56e-02 \\ 
\end{tabular}
\label{tab:circular_convection_L1}
\end{subtable}
\par\bigskip
\begin{subtable}{\textwidth}
\centering
\caption{Error in $\left\Arrowvert \cdot \right\Arrowvert_2$.}
\begin{tabular}{c c c c}
\textbf{Level} & \textbf{MC} & \textbf{LP} & \textbf{MU} \\ \hline
5 & 1.46e-01 & 9.76e-02 & 1.61e-01 \\ 
6 & 8.15e-02 & 6.64e-02 & 8.52e-02 \\ 
7 & 5.81e-02 & 4.59e-02 & 5.60e-02 \\ 
\end{tabular}
\label{tab:circular_convection_L2}
\end{subtable}
\label{tab:circular_convection_errors}
\end{table}